\documentclass[reqno,11pt]{amsart}
\usepackage{amsmath,amssymb} 
\usepackage{amsthm, amsfonts, amsgen, stmaryrd,enumerate}
\usepackage[english]{babel}
\usepackage{fullpage}
\usepackage[centertags]{amsmath}
\usepackage{indentfirst}
\usepackage{fancyhdr}
\usepackage[dvips]{graphicx}
\usepackage{psfrag}
\numberwithin{equation}{section}
\usepackage[latin1]{inputenc}

\usepackage{color}

\title[Hodge dualities on  quantum spheres]{Hodge duality operators on left covariant exterior algebras over  two and three  dimensional quantum  spheres}
\date{5 June 2012}
\author{Alessandro Zampini}
\address{ Mathematisches Institut der L.M.U., Theresienstra\ss e 39,  D-80333 M\"unchen, Germany.}
 \email{zampini@math.lmu.de}

\newtheorem{theo}{Theorem}[section]

\newtheorem{rema}[theo]{Remark}

\newcommand{\ii}{\mathrm{i}}
\newcommand{\nn}{\nonumber}

\newcommand{\dd}{{\rm d}}
\newcommand{\ca}{\mathcal{A}}

\newcommand{\ch}{\mathcal{H}}
\newcommand{\cl}{\mathcal{L}}

\newcommand{\cq}{\mathcal{Q}}

\newcommand{\cu}{\mathcal{U}}        
\newcommand{\SU}{\mathrm{SU}_q(2)}  
\newcommand{\ASU}{\ca(\mathrm{SU}_q(2))}  
\newcommand{\sq}{\mathrm{S}^2_{q}}  
\newcommand{\Asq}{\ca(\mathrm{S}^2_{q})}  
\newcommand{\su}{\cu_q(\mathfrak{su}(2))}  

\newcommand{\eps}{\varepsilon}      
\newcommand{\hs}[2]{\left\langle #1,#2\right\rangle}  
\newcommand{\hp}[2]{\left\{ #1,#2\right\}} 

\newcommand{\lt}{{\triangleright}}    
\newcommand{\rt}{{\triangleleft}}
\newcommand{\IC}{{\mathbb C}} 
\newcommand{\IR}{{\mathbb R}} 
\newcommand{\IN}{{\mathbb N}} 
\newcommand{\IZ}{{\mathbb Z}} 
\DeclareMathOperator{\id}{id}       
\DeclareMathOperator{\U}{U}       

\newcommand{\figureheight}{8cm}
\newcommand{\putfig}[2]{\begin{figure}[htp]
        \special{isoscale c:/itex/texfig/#1.wmf, \the\hsize \figureheight}
        \vspace{\figureheight}
        \caption{#2}\label{fig:#1}
        \end{figure}}
\newcommand{\pictureheight}{4cm}
\newcommand{\putpicture}[2]{\begin{figure}[htp]
        \special{isoscale c:/itex/texfig/#1.wmf, \the\hsize \pictureheight}
        \vspace{\pictureheight}
        \caption{#2}\label{fig:#1}
        \end{figure}}

\newcommand{\beqa}{\begin{eqnarray}}
\newcommand{\eeqa}{\end{eqnarray}}
\newcommand{\beq}{\begin{equation}}
\newcommand{\eeq}{\end{equation}}

\newcommand{\complex}{{\mathbb C}} 

\newcommand{\oz}{\omega_{z}}
\newcommand{\op}{\omega_{+}}
\newcommand{\om}{\omega_{-}}
\newcommand{\ot}{\otimes}

\newcommand{\gm}{\mathrm{g}}
\newcommand{\cx}{\mathcal{X}}
\newcommand{\ck}{\mathcal{K}}

\begin{document}

\thispagestyle{empty}

\begin{abstract}
Using  non canonical braidings, we first introduce a notion of symmetric tensors and corresponding Hodge operators on a class of left-covariant 3d differential calculi over $\SU$, then we induce Hodge operators on the left covariant 2d exterior algebra over the Podle\`s quantum sphere.    
\end{abstract}


\maketitle
\tableofcontents

\section{Introduction}

Following the formalism developed by Woronowicz \cite{wor89}, various  aspects of the differential geometry induced on a  large class of quantum groups $\ch$ equipped with  suitable bicovariant differential first order calculi $(\dd, \Gamma)$  have been intensively studied. 

Among those, the general problem of defining meaningful symmetric tensors and Hodge duality operators acting on the higher order differential calculi $\Gamma_{\sigma^{\pm}}$ constructed via the canonical braiding $\sigma$ on $\Gamma^{\otimes2}=\Gamma\otimes_{\ch}\Gamma$ (and its inverse  $\sigma^{-1}=\sigma^{-}$,  the braiding being no longer the classical flip) for the calculus has been considered in detail in \cite{hec99,hec00,hec03}. 
Given finite $N$-dimensional calculi such that $\dim\,\Gamma^k_{\sigma^{\pm}}\,=\,\dim\,\Gamma^{N-k}_{\sigma^{\pm}}$ for each space of exterior $k$- and $(N-k)$-forms, the properties of the canonical braiding and of its corresponding antisymmetriser operators $A^{(k)}_{\sigma^{\pm}}$ on $\Gamma_{\sigma^{\pm}}^k$  allow to define a rank two symmetric tensor over $\Gamma^{\otimes2}$ and Hodge operators $\star_{\sigma^{\pm}}:\Gamma^{k}_{\sigma^{\pm}}\to\Gamma^{N-k}_{\sigma^{\pm}}$  satisfying $\star_{\sigma^{\pm}}\star_{\sigma^{\mp}}=1$. 

It is crucial in this formulation that the first order differential calculus is bicovariant, a condition which is sufficient to have a canonical braiding \cite{ks}. The properties of the canonical braiding  also core  the formalism developed  in \cite{maj-ri-99, maj-ri-02} to introduce the concept of Riemannian quantum group and braided Killing form as a part of a more general formulation aimed to describe framings and coframings over quantum groups as a gauge theory of quantum differential forms. 

A  somehow reversed strategy for the specific example of the quantum ${\rm SU(2)}$ group equipped with the bicovariant $4D_+$ calculus gives in  \cite{ale11} results which largely  agree to (and slightly generalize)  those presented in the former approach. A Hodge operator is there meant as a bijection $\star_{\sigma^{\pm}}:\Gamma^{k}_{\sigma^{\pm}}\to\Gamma^{N-k}_{\sigma^{\pm}}$ whose square has, for a  suitably defined  class of symmetric tensors, the same degeneracy of the antisymmetrisers of the calculus.  

The bicovariance of the calculus seems to play in this approach no explicit role, and this suggests  that it is  possible to study the problem of defining Hodge duality operators even on exterior algebras built over left-covariant calculi on quantum groups, provided they have a consistent --  although not canonical -- braiding.  Presenting the first  results obtained  in this direction is one of the aim of the present paper. 

Using the classification of \cite{hec2001}, we equip the quantum group $\SU$ with a set (that we call $\ck$)  of left covariant first order  3d calculi having a non canonical braiding $\sigma$ and study the exterior algebras associated to the corresponding  antisymmetriser operators. We introduce then Hodge operators acting on such exterior algebras: by 
 \emph{consistent}  we mean that their squares present  the same degeneracy of the antisymmetriser operators. 
Such Hodge operators exist for a class of properly defined symmetric tensors $\gm$ acting on $\Gamma^{\otimes2}$, with a notion of symmetry which does not necessarily coincide with the standard $\gm\,\circ\,\sigma\,=\,\gm$. 
The class $\ck$ of calculi we consider contains the most famous Woronowicz' 3d calculus \cite{wor87}: 
this paper then deepens  the analysis on scalar products and duality operators presented in \cite{funingeo}, and enlarges its results.

That the formalism we develop is consistent for any calculus  in $\ck$ poses the further question this paper aims to analyze. Is it possible to introduce a notion eventually selecting a proper and interesting subclass of elements in $\ck$, i.e. of calculi on $\SU$ among those that are being considered? We propose this condition to be the requirement that the previously introduced notion of symmetry for a tensor $\gm$ does coincide with the standard one, namely that $\gm\,\circ\,\sigma\,=\,\gm$. 
This condition actually selects a subset $\tilde{\ck}\,\subset\,\ck$ of calculi on $\SU$, a subset  we can describe following a further interesting characterization. 

It is well known that the Podle\`s standard sphere $\sq$ is the quantum homogeneous space defined by a $\U(1)$-coaction over $\SU$; for a set $\ck_{\pi}\,\subset\,\ck$ of so called projectable calculi on $\SU$ this topological Hopf fibration acquires compatible differential structures. The restriction to $\sq$ of such projectable calculi gives in particular different isomorphic realizations of the unique 2d left covariant differential calculus introduced by Podle\`s himself \cite{po92}. 
Considering any of these realizations of the exterior algebra  $\Gamma(\sq)$ in terms of the frame bundle approach \cite{maj-ri-02}  allows to describe how the restriction of the above Hodge operators (on $\SU$) meaningfully introduces a class (since they correspond to elements in $\ck_{\pi}$) of different bijections as maps $\check{S}\,:\,\Gamma^{k}(\sq)\,\to\,\Gamma^{2-k}(\sq)$. It turns out that the operators $\check{S}$ have the expected degeneracy (i.e. the degeneracy of a Hodge duality on a classical 2 dimensional exterior algebra) only if they are induced by the Hodge operators on $\SU$ equipped with the calculi $\tilde{\ck}\,\subset\,\ck$. 

It means (this being the last result that this paper presents) that the formalism developed here allows to induce Hodge operators acting on the left invariant 2d exterior algebra $\Gamma(\sq)$ starting from a suitable formulation of symmetric tensors and Hodge operators on the quantum group $\SU$.

The paper is organised as follows. Section \ref{se:gese} describes  the geometrical setting of the analysis, namely those aspects of differential calculi and exterior algebras over classical and quantum group we shall use, and present the class $\ck$ of differential calculi over the quantum  ${\rm SU(2)}$ we shall consider.  Starting from a tensor whose coefficients give  contraction map, section \ref{hosc:s} present families of scalar products and corresponding dual Hodge operators. The example of the Woronowicz' calculus is used as a guide, the results are then extended to the whole class $\ck$ of calculi. In section 
\ref{shosf} we finally study Hodge operators for the quantum sphere $\sq$.

\section{The geometrical setting}
\label{se:gese}

\subsection{A classical setting}
\label{ove:ss}

Consider a $N$-dimensional connected Lie group $G$ given as the real form of a complex connected Lie group. Its group manifold is parallelizable: the space of 1-forms $\Omega^1(G)$ is a free bicovariant $N$-dimensional  $\ca(G)$-bimodule on the basis of left (right) invariant $\{\phi^{a}\}$ ($\{\eta^{a}\}$) 1-forms. 
The associated first order differential calculus is given by $(\dd, \Omega^1(G)$) with the exterior differential given by $\dd h=(L_{a}h)\omega^{a}=(R_{a}h)\eta^a$ in terms of the action of the dual  left (right) invariant derivations $L_{a}$ ($\{R_{a}\}$) on $h\in\,\ca(G)$. 

The standard flip given\footnote{We denote $\Omega^{\otimes k}=\Omega^{1}(G)\otimes_{\ca(G)}\cdots\otimes_{\ca(G)}\Omega^1(G)$ and drop the overall obvious  dependence on $G$.}  on a basis by $\tau:\omega^{a}\otimes\omega^b\,\mapsto\,\omega^b\otimes\omega^a$ is a braiding on $\Omega^{\otimes 2} $; the corresponding  antisymmetriser operators $A^{(k)}\,:\,\Omega^{\otimes k}\,\to\,\Omega^{\otimes k}$  ($k\,\in\,\IN$) give the exterior algebra $\Omega^{\wedge}=(\oplus_{k=1}^{N}\,\Omega^{k}, \wedge)$ as 
$\Omega^{\otimes k}\,\supset\,\Omega^k\,=\,
(\mathrm{Range}\,A^{(k)})\,\simeq\,\Omega^{\otimes k}/\ker\,A^{(k)}$, with 
\beq
\label{dAom}
\omega^{a_{1}}\wedge\ldots\wedge\omega^{a_{k}}\;=\;A^{(k)}(\omega^{a_{1}}\otimes\ldots\otimes\omega^{a_{k}})\;=\;\sum_{\pi\in\,S_{k}}(-1)^{\pi}\omega^{\pi(a_{1})}\otimes\ldots\otimes\omega^{\pi(a_{k})}
\eeq
 where  $S_{k}$ is the set of permutations of $k$ elements. 
The differential calculus $(\Omega^{\wedge}, \dd)$ is given by equipping the exterior algebra with the unique consistent graded derivative operator $\dd:\Omega^{ k}\,\to\,\Omega^{k+1}$  satisfying $\dd^2=0$ and a graded Leibniz rule. Every $\Omega^{k}$ is a bicovariant free $\ca(G)$-bimodule with $\dim\,\Omega^{ k}=N!/(k!(N-k)!)$ and $\Omega^{ k}=\emptyset$ for $k>N$. 
The antisymmetrisers have a completely degenerate spectral decomposition,
\beq
\label{spclA}
A^{(k)}(\omega^{a_{1}}\wedge\ldots\wedge\omega^{a_{k}})=k!(\omega^{a_1}\wedge\ldots\wedge\omega^{a_{k}}).
\eeq
We consider a non degenerate tensor $\gm\,:\,\Omega^1\times\Omega^1\,\to\,\ca(G)$, whose components we use to set an  $\ca(G)$-bimodule contraction $\gm\,:\,\Omega^{\otimes k}\times\Omega^{\otimes (k+k^{\prime})}\,\to\,\Omega^{k^{\prime}}$ given on a  basis by
\beq
\gm(\omega^{a_{1}}\otimes\ldots\otimes\omega^{a_{k}},\omega^{b_{1}}\otimes\ldots\otimes\omega^{b_{k+k^{\prime}}})\,=\,
\left\{\Pi_{j=1,\ldots, k}\,\gm(\omega^{a_j},\omega^{b_{j}})\right\}\,\omega^{b_{k+1}}\otimes\ldots\otimes\omega^{b_{k^{\prime}}};
\label{clac}
\eeq
the properties of the antisymmetriser operators allow then to prove that the position (see \eqref{dAom})
\beq
\label{clacon}
\gm(\omega^{a_{1}}\wedge\ldots\wedge\omega^{a_{k}},\omega^{b_1}\wedge\ldots\wedge\omega^{b_{k+k^{\prime}}})\,=\,g(A^{(k)}(
\omega^{a_{1}}\otimes\ldots\otimes\omega^{a_{k}}),A^{(k+k^{\prime})}(\omega^{b_{1}}\otimes\ldots\otimes\omega^{b_{k+k^{\prime}}}))
\eeq
consistently generalizes the contraction map \eqref{clac} to $\gm:\Omega^k\times\Omega^{k+k^{\prime}}\,\to\,\Omega^{k^{\prime}}$. 
If  $\mu\,=\,m\,\theta\,=\,\mu^*$ is a volume form, with $\theta=\omega^1\wedge\ldots\wedge\omega^N$ the top form corresponding to an ordering of the basis elements $\{\omega^a\}$ and  $m\,\in\,\IC$), we define the operator $S:\Omega^k\to\Omega^{N-k}$ by 
\beq
\label{fueq}
S(\phi)\,=\,\frac{1}{k!}\,\gm(\phi,\mu),
\eeq
on any $k$-form $\phi$.  The following equivalence holds ($\phi,\phi^{\prime}\,\in\,\Omega^1$)
\beq
\label{claeq}
\gm(\phi,\phi^{\prime})\,=\,\gm(\phi^{\prime},\phi)\qquad\Leftrightarrow\qquad
S^2(\phi)=(-1)^{N-1}\{S^2(1)\}\phi.
\eeq
The tensor $\gm$ is symmetric if and only if the action of   the restiction of $S^2$ on $\Omega^1$ is a constant depending on the volume\footnote{It is true that the factor can be  arbitrary, and that $\gm$ is symmetric if and only if $S^2(\phi^a)=\zeta\,\phi^a$ with $0\,\neq\,\zeta\,\in\,\IC$.  The choice in \eqref{claeq} will give the possibility of the usual overall normalization.};  given a symmetric $\gm$ one has also that $S^{2}(\phi)=(-1)^{k(N-k)}\{S^2(1)\}\phi$ with $\phi\,\in\,\Omega^k$. But such an operator $S$ is not (yet) an Hodge operator: it has to be real, and  the reality condition comes as  the equivalence
\beq
\label{reacla}
\gm(\phi,\phi^{\prime})^*\,=\,\gm(\phi^{\prime*},\phi^{*})\qquad\Leftrightarrow\qquad S(\phi^{*})\,=\,(S(\phi))^*.
\eeq
The compatibility of the action of $S$ with the hermitian conjugation on $\Omega^1$ and $\Omega^{N-1}$  turns out to be  sufficient to have $[S,^*]=0$ on the whole exterior algebra $\Omega^{\wedge}$. Such a symmetric and real operator $S$ is then recovered as the Hodge operator corresponding to the (inverse) of the (metric) tensor $\gm$ on the group manifold. The choice  $S^2(1)\,=\,{\rm sgn}(\gm)$ fixes the modulus of the scale parameter $m$ of the volume so to have 
$$
S^2(\phi)\,=\,{\rm sgn}(\gm)(-1)^{k(N-k)}\phi
$$ 
for any $\phi\,\in\,\Omega^k$.

Hodge operators can be introduced also on homogeneous spaces. 
Let $K\,\subset\,G$ be a compact Lie subgroup of $G$. The quotient of its right action ${\rm r}_{k}(g)\,=\,g\,k$ for $k\,\in\,K$ and $g\,\in\,G$ gives a principal fibration $\pi\,:\,G\,\to\,G/K$. A homogeneous space is not necessarily parallelizable: the exterior algebra $\Omega(G/K)\,\subset\Omega(G)$ is given by horizontal and right $K$-invariant forms on $G$,
\beq
\label{hoscl}
\Omega(G/K)\,=\,\{\psi\,\in\,\Omega(G)\,:\,i_{X_{V}}\psi\,=\,0\,;\, {\rm r}_{k}^*(\psi)\,=\,\psi\},
\eeq
with $X_{V}$ the vertical fields of the fibrations (i.e. the infinitesimal (left-invariant) generators of the right $K$ action on $G$),  and ${\rm r}_{k}^*$ the natural pull-back action to $\Omega(G)$.
The $\Omega^{s}(G/K)$ sets (with $0\,<\,s\,<\,N^{\prime}$) are no longer  free $\ca(G/K)$-bimodules;   
the dimension of the exterior algebra $\Omega(G/K)$ is  given as the highest integer $N^{\prime}$ so that $\Omega^{N^{\prime}+1}(G/K)\,=\,\emptyset$, and coincides with $N^{\prime}\,=\,\dim G\,-\,\dim K$. 
The set $\Omega^{N^{\prime}}(G/K)$ is indeed a free 1-dimensional $\ca(G/K)$ bimodule with a basis element given by $\check{\theta}\,=\,i_{X_{V^1}}\cdots i_{X_{V^{\dim K}}}\theta$ for a basis  $X_{V^{a}}$ of the Lie algebra of vertical vector fields of the fibration. We have then a consistent (up to scalars) left invariant volume form $\check{\mu}\,=\,\check{\mu}^*$ on the homogeneous space: if  we consider right $K$-invariant metric tensors $\gm$ on $G$ whose restriction to the homogeneus space $G/K$ is non degenerate, then  the map $\check{S}\,:\,\Omega^{j}(G/K)\,\to\,\Omega^{N^{\prime}-j}(G/K)$ given by
\beq
\check{S}(\psi)\,=\,\frac{1}{j!}\,\gm(\psi,\check{\mu})
\label{Ssf}
\eeq
is a well-defined bijection, satisfying the relation $\check{S}^2(\psi)\,=\,{\rm sgn }(\gm(\check{\mu},\check{\mu}))(-1)^{s(N^{\prime}-s)}\psi$ for any $\psi\,\in\,\Omega^{k}(G/K)$ after a natural normalisation.

Both the Hodge operators above can be formulated following a different path. 
Starting from a non degenerate  tensor $\gm$, a sesquilinear map $\hs{~}{~}_{G}:\Omega^k\times\Omega^k\,\to\,\ca(G)$ can be defined by
 \beq
 \label{stancooggi}
 \hs{\phi}{\phi^{\prime}}_{G}\,=\,\frac{1}{k!}\,\gm(\phi^*,\phi^{\prime});
 \eeq
   the equation
 \beq
 \label{claTop}
\phi^{*}\wedge T(\phi^{\prime})\,=\,(\hs{\phi}{\phi^{\prime}}_{G})\mu
 \eeq
uniquely defines a bijective $T\,:\,\Omega^{k}\,\to\,\Omega^{N-k}$ with $T(1)=\mu, \;T(\mu)=m$. It is immediate to check the equivalence $T(\phi)\,=\,S(\phi)$ on any $\phi\,\in\,\Omega^k$, which comes  from 
\beq
\gm(\phi,\phi^{\prime})\,\mu\,=\,\phi^{*}\wedge\gm(\phi^{\prime},\mu).
\label{nonso}
\eeq
for any pair $\phi,\phi^{\prime}\,\in\,\Omega^k$.
It 
is analogously immediate to see that the 
restriction
$$
\hs{\psi}{\psi^{\prime}}_{G/K}\,=\,\hs{\psi}{\psi^{\prime}}_{G}
$$
of the  sesquilinear allows to consistently set 
\beq
\label{hoho}
\psi^{*}\wedge\check{T}(\psi^{\prime})\,=\,(\hs{\psi}{\psi^{\prime}}_{G/K})\check{\mu}
\eeq
as a definition for the operator $\check{T}\,:\,\Omega^{j}(G/K)\,\to\,\Omega^{N^{\prime}-j}(G/K)$.  One has clearly $\check{T}=\check{S}$  as the Hodge operators on $\Omega(G/K)$ corresponding to  projecting the  right $K$-invariant (inverse) metric  tensor $\gm$ onto the homogeneous space.


\subsection{A quantum setting: left covariant differential calculi over quantum groups}
\label{ss:qs}

Consider  $\ch$ to be the unital $*$-Hopf algebra  $\ch=(\ch,\Delta,\varepsilon,S)$ over $\IC$, with  $\Gamma$ an $\ch$-bimodule. The pair $(\Gamma, \dd)$ is a (first order) differential calculus over $\ch$  provided the  linear map $\dd:\ch\to\Gamma$ satisfies the Leibniz rule, $\dd (h\, h^{\prime}) = (\dd  h) h^{\prime}  + h\, \dd h^{\prime} $ for $h,h^{\prime}\in \ch$, and 
$\Gamma$ is generated by $\dd(\ch)$ as a $\ch$-bimodule. It is called  a $*$-calculus provided  there is an anti-linear involution $*:\Gamma\to\Gamma$ such that $(h_{1}(\dd h)h_{2})^*=h_{2}^{*}(\dd(h^*))h_{1}^{*}$ for any $h,h_{1},h_{2}\in \ch$. 

A first order differential calculus is said left covariant provided  a left coaction $\Delta_{L}^{(1)}:\Gamma\to\ch\otimes\Gamma$ exists, such that $\Delta_{L}^{(1)}(\dd h)=(1\otimes \dd)\Delta(h)$ and $\Delta_{L}^{(1)}(h_1\,\alpha\,h_2)=\Delta(h_{1})\Delta_{L}^{(1)}(\alpha)\Delta(h_2)$ for any $h,h_1,h_2\in\,\ch$ and $\alpha\in\,\Gamma$. The set $\Gamma$ turns out to be a free left covariant $\ch$-bimodule, with a free basis $\Gamma_{\rm L}$ of left invariant one forms,   namely the elements $\omega_{a}\in \Gamma$ such that $\Delta^{(1)}_{L}(\omega_{a})=1\otimes\omega_{a}.$ Its dimension is called the dimension of the first order calculus. The map $\mathfrak{R}\,:\,\ch\,\to\,\Gamma_{\rm L}$ given by
\beq
\label{mappafr}
\mathfrak{R}(h)\,=\,S(h_{(1)})\,\dd h_{(2)}
\eeq
allows to characterise left covariant first order  differential calculi:  they correspond to the choice of a right ideal 
$\mathcal{Q} \subset \ker\eps$ with 
\beq
\label{cri}
\cq\,=\,\{h\,\in\,\ker\,\eps\,:\,\mathfrak{R}(h)=0\}; 
\eeq
there is a left $\ch$-modules isomorphism given by 
$\Gamma\simeq\ch\otimes(\ker\eps/\cq)$, and a complex vector space isomorphism $\Gamma_{\rm L}\simeq\ker\eps/\cq$.

The tangent space of
the calculus is  the complex vector space of elements  out of $\ch^{\prime}$ -- the dual space $\ch^{\prime}$ of functionals on $\ch$ --  defined by 
$\mathcal{X}_{\mathcal{Q}}:=
\{X\in \ch^{\prime} ~:~ X(1)=0,\,X(Q)=0, \,\, \forall \, Q\in\mathcal{Q}\}. $ One has that $(\Gamma,\dd)$ is a $*$ calculus if and only if its quantum tangent space is $*$-invariant.
There exists a unique bilinear form 
\beq
\hp{~}{~}:\mathcal{X}_{\cq}\times\Gamma, \qquad\hp{X}{x \dd y}:=\eps(x)X(y) ,
\label{dptg}
\eeq
giving a non-degenerate dual pairing between the vector spaces $\mathcal{X}_{\cq}$ and $\Gamma_{\rm L}$. The dual space $\ch^{\prime}$ has natural left and right (mutually commuting) actions  on $\ch$:  
\beq
X\triangleright h:=h_{(1)}X(h_{(2)}),\qquad
 h\triangleleft X:=X(h_{(1)})h_{(2)}.
\label{deflr}
\eeq
If the  vector space  $\mathcal{X}_{\cq}$ is finite dimensional,  its elements belong to the dual Hopf algebra 
$ \ch^{\prime} \supset \ch^{o} = (\ch^{o}, \Delta_{\ch^o}, \eps_{\ch^o},S_{\ch^o})$, defined as the largest Hopf $*$-subalgebra contained in $\ch^{\prime}$. 
In such a case  the 
$*$-structures are compatible with both actions,
$$
X \lt h^* = ((S(X))^* \lt h)^*,\qquad
h^* \rt  X = (h \rt  (S(X))^*)^*,
$$
for any $ X \in \ch^{o}, \ h \in \ch$ 
and the exterior derivative can be written as:
\beq
\dd h := \sum_a
~(X_{a} \triangleright h) ~\omega_{a} =\sum_{a} \omega_{a}(-S^{-1}(X_{a}))\lt h, 
\label{ded}
\eeq
where $\hp{X_{a}}{\omega_{b}}=\delta_{ab}$.   
The twisted Leibniz rule of derivations of the basis elements $X_{a}$ is dictated by their coproduct:  
\beq
\Delta_{\ch^o}(X_{a})=1\otimes X_{a}+
\sum\nolimits_b X_{b}\otimes f_{ba},
\label{cpuh}
\eeq
where the $f_{ab} \in \ch^{o}$ consitute an algebra representation of $\ch$, also  controlling the $\ch$-bimodule structure of $\Omega^{1}(\ch)$:
\begin{align}\label{bi-struct}
\omega_{a} h = \sum\nolimits_b (f_{ab}\triangleright h)\omega_{b} , \qquad 
h \omega_{a} = \sum\nolimits_b \omega_{b} \left( (S^{-1}(f_{ab}) )\triangleright h \right) , 
\qquad \mathrm{for} \quad h \in \ch. 
\end{align}
In order to build an exterior algebra over the FODC $(\dd, \Gamma)$, consider $\Gamma^{\otimes k}$ as the k-fold tensor product $\Gamma\otimes_{\ch}\cdots\otimes_{\ch}\Gamma$ (with $\Gamma^{0}=\ch$) and $\Gamma^{\otimes}=\oplus_{k=0}^{\infty}\Gamma^{\otimes k}$, which is an algebra with multiplication $\otimes_{\ch}$. From the map 
\beq
\mathcal{S}\,:\,\ch\to\Gamma_{\rm L}^{\otimes 2},\qquad\qquad x\,\mapsto\,\sum\mathfrak{R}(x_{(1)})\otimes\mathfrak{R}(x_{(2)}),
\label{defSS}
\eeq
let $\mathcal{S}_{\mathcal{Q}}\,\subset\,\Gamma^{\otimes2}$ be the 2-sided ideal in $\Gamma^{\otimes}$ generated by the  range of its restriction to $x\,\in\,\mathcal{Q}$; the quotient $\Gamma_{u}^{k}=\Gamma^{\otimes k}/(\mathcal{S}_{\mathcal{Q}}\cap\Gamma^{\otimes k})$ is a well defined $\ch$-bimodule.

This exterior algebra turns out to be a differential calculus over $\ch$ once the exterior derivative $\dd$ is extended as a graded derivation with $\dd^2=0$,  satisfying a graded Leibniz rule (that is $\dd(\omega\wedge\omega^{\prime})=(\dd\omega)\wedge\omega^{\prime}+(-1)^{m}\omega\wedge\dd\omega^{\prime}$ for any $\omega\in\,\Gamma^{m}$). The quotient $\Gamma_{u}$ also inherits the natural extension of the left coaction of $\ch$, which is compatible with the action of the operator $\dd$, so to have a left covariant differential calculus $(\dd, \Gamma_{u})$ over the FODC which is universal: any other left covariant differential calculus $(\dd, \Gamma)$ over $\ch$ with $\Gamma^{1}=\Gamma$ is    a suitable quotient of the universal one.

Given the left covariant bimodule $\Gamma$ over $\ch$, an invertible linear mapping $\sigma:\Gamma\otimes_{\ch}\Gamma\,\to\,\Gamma\otimes_{\ch}\Gamma$ is called a braiding for $\Gamma$ provided $\sigma$ is a $\ch$-bimodule homomorphism which commutes with the left coaction on $\Gamma$ and satisfies the braid equation 
\beq
(1\otimes\sigma)\circ(\sigma\otimes 1)\circ(1\otimes \sigma)=(\sigma\otimes 1)\circ(1\otimes \sigma)\circ(\sigma\otimes 1)
\label{br3}
\eeq
on $\Gamma^{\otimes3}$. The next natural requirement is that $(1-\sigma)(\mathcal{S}_{\mathcal{Q}}\cap\Gamma^{\otimes2})=0$. Such a braiding neither needs to exist nor it is unique for a given left covariant differential calculus over $\ch$: this is the main difference with bicovariant differential calculi, which present a  canonical braiding.  If a braiding does exists, then the corresponding antisymmetriser operators $A^{(k)}:\Gamma^{\otimes k}\to\Gamma^{\otimes k}$ are well defined and their ranges give the differential calculus $(\dd, \Gamma_{\sigma})$ since $\ker\,A^{(k)}\supset\mathcal{S}_{\mathcal{Q}}$ is a 2-sided graded ideal in $\Gamma^{\otimes k}$.

\subsection{A class of left covariant  differential calculi over the quantum SU(2)}
\label{ss:sesu}

As quantum group $\SU$ we consider the compact real form of the quantum group ${\rm SL}_{q}(2)$ and, following \cite{worles},  we formulate it as the  polynomial  unital  $*$-algebra $\ASU=(\SU,\Delta,S,\eps)$ generated by elements $a$ and $c$ which we write using the matrix notation
\beq
\label{Us}
 u = 
\left(
\begin{array}{cc} a & -qc^* \\ c & a^*
\end{array}\right).
\eeq
The Hopf algebra structure can then be expressed as
$$
uu^*=u^*u=1,\quad \,\Delta\, u = u \otimes u ,  \quad S(u) = u^* , \quad \eps(u) = 1
$$
 with the deformation parameter
$q\in\IR$. 


In order to describe the quantum tangent spaces of the calculi that will be later  introduced,  we consider  the set of functionals given by the unital Hopf $*$-algebra $\widetilde{\ASU}$ over $\complex$, satisfying the inclusions $\ASU^{o}\supset\widetilde{\ASU}\supset\su$ with $\ASU^{o}$ the Hopf dual $*$-algebra  and $\su$ the universal envelopping algebra of $\ASU$.
As an algebra is $\widetilde{\ASU}$ generated by  the five elements $\{K^{\pm 1},E,F,\eps_{-}\}$, with $K K^{-1}=1$ fullfilling the
relations\footnote{We shall also denote $K^{+}=K, \,K^{-}=K^{-1}$. It is clear that $\widetilde{\ASU}$ is generated by the universal envelopping $\su$ algebra together with the $\ASU$-character $\eps_{-}$ acting as $\eps_{-}(a)=\eps_{-}(a^*)=-1; \;\eps_{-}(c)=\eps_-(c^*)=0$. }: 
\begin{align}
&\eps_{-}K^{\pm}=K^{\pm}\eps_{-},\qquad\eps_{-}\eps_{-}=1, \nonumber \\
&\qquad \eps_{-}E=E\eps_{-},\qquad\eps_{-}F=F\eps_{-},\nonumber \\
 &K^{\pm}E=q^{\pm}EK^{\pm}, \qquad 
K^{\pm}F=q^{\mp}FK^{\pm}, \nonumber \\
& \qquad  
[E,F] =\frac{K^{2}-K^{-2}}{q-q^{-1}}.
\label{relsu}
\end{align} 
The $*$-structure is
$K^*=K, \,  E^*=F ,\,\eps_{-}^*=\eps_{-}$,
while the Hopf algebra structures are 
\begin{align*}
&\Delta(K^{\pm}) =K^{\pm}\otimes K^{\pm}, \qquad
\Delta(E) =E\otimes K+K^{-1}\otimes E,  \\ &\quad 
\Delta(F)
=F\otimes K+K^{-1}\otimes F,
\qquad \Delta(\eps_{-})=\eps_{-}\otimes\eps_{-},
\\ &\quad\quad S(K) =K^{-1}, \qquad
S(E) =-qE, \qquad 
S(F) =-q^{-1}F , \qquad S(\eps_{-})=\eps_{-}, \\ 
&\quad\quad\quad\varepsilon(K)=\varepsilon(\eps_{-})=1, \quad \varepsilon(E)=\varepsilon(F)=0.
\end{align*}
The only non zero terms of its action on  $\ASU$ is given on the generators by
\begin{align}
&K^{\pm}(a)=q^{\mp 1/2}, \qquad  K^{\pm}(a^*)=q^{\pm 1/2},  \qquad E(c)=1, \qquad F(c^*)=-q^{-1}, \nn \\
&\qquad\eps_{-}(a)=\eps_{-}(a^*)=-1.\label{ndp}
\end{align}

Given the algebra $\ca(\U(1)):=\IC[z,z^*] \big/ \!\!<zz^* -1>$,  the map  
\beq \label{qprp}
\pi: \ASU \, \to\,\ca(\U(1)),\quad\quad\pi(a)\,=\,z, \quad\pi(a^*)\,=\,z^*,\quad\pi(c)\,=\,\pi(c^*)\,=0
\eeq 
is a surjective Hopf $*$-algebra homomorphism, so that  $\U(1)$
is a quantum subgroup of $\SU$ with right coaction:
\beq 
\delta_{R}:= (\id\otimes\pi) \circ \Delta \, : \, \ASU \,\to\,\ASU \otimes
\ca(\U(1)) . \label{cancoa} 
\eeq 
This right coaction gives a decomposition 
\beq
\label{clnm}
\ASU=\oplus_{n\in\,\IZ}\cl_{n}, \qquad\qquad\cl_{n}:=\{x\in\,\ASU\;:\;\delta_{R}(x)=x\otimes z^{-n}\},
\eeq
with $\Asq=\cl_{0}$ the algebra of the standard Podle\'s sphere, each $\Asq$-bimodule $\cl_{n}$ giving the set of (charge $n$)  $\U(1)$-coequivariant maps for the topological quantum principal bundle $\Asq\hookrightarrow\ASU$.  

\bigskip

From \cite{hec2001} we know a classification of left covariant differential calculi over the  quantum group ${\rm SL}_{q}(2)$. Among them we select those calculi which are compatible  with the reality structure (the anti-hermitian involution) giving $\SU$, and  which present a consistent braiding. This means that we equip $\SU$ with the left covariant differential calculi satisfying the following properties: 
~\\
\begin{itemize}
\item $\Gamma$ is a left covariant bimodule;
~\\
\item a basis of $\Gamma_{\rm L}$ is given by $\mathfrak{R}(c), \mathfrak{R}(c^*), \mathfrak{R}(a-a^*)$;
~\\ 
\item for the corresponding universal differential calculus it is $\dim \Gamma_{u}^{\wedge 2}\geq 3$; 
~\\
\item $\Gamma$ is Hopf-invariant, i.e. for the corresponding right ideal $\cq\subset\ker\eps$ the equality $\varphi(\cq)=\cq$ for any Hopf algebra automorphism $\varphi$ on $\ASU$ holds.
~\\
\item On $\Gamma^{\otimes2}$ a consistent braiding exists.
\end{itemize}
~\\
We are then left with seven (up to isomorphisms) such calculi, and we denote this class by $\ck$. We present these calculi giving  \eqref{cri} the generators of  the right ideals $\cq_{j}\subset\ker\eps$ (with $j=1,\ldots,7$) they are characterised by, together with a  basis of their quantum tangent spaces:
~\\
\begin{enumerate}
\item $\cx_{\cq_1}=\{X_{z}\,=\,\frac{K^{-2}-1}{q-1}, \quad X_{+}\,=\,q^{-1/2}EK^{-1}, \quad X_{-}\,=\,q^{1/2}FK^{-1}\}$,
\begin{align}\cq_{1}\,=\{a+qa^*-(1+q), c^2, c^{*2}, cc^*, (a-q)c, (a-q)c^*\};\label{q1} \\ \nn \end{align}
\item $\cx_{\cq_2}=\{X_{z}\,=\,\frac{\eps_{-}K^{2}-1}{q+1}, \quad X_{+}\,=\,-q^{-1/2}\eps_{-}EK^{-1}, \quad X_{-}\,=\,-q^{1/2}\eps_{-}FK^{-1}\}$,  
\begin{align}
 \cq_{2}\,=\{a-qa^*-(1-q), c^2, c^{*2}, cc^*, (a+q)c, (a+q)c^*\};
 \label{q2} 
 \\ \nn
 \end{align}
\item $\cx_{\cq_3}=\{X_{z}\,=\,\frac{K^{-4}-1}{q^2-1}, \quad X_{+}\,=\,q^{-3/2}EK^{-3}, \quad X_{-}\,=\,q^{3/2}FK^{-3}\}$, 
\begin{align}
\cq_{3}\,=\{a+q^2a^*-(1+q^2), c^2, c^{*2}, cc^*, (a-q^2)c, (a-q^2)c^*\};
\label{q3} 
\\ \nn
\end{align}
\item $\cx_{\cq_{4}}=\{X_{z}\,=\,\frac{K^{2}-1}{q^{-1}-1}, \quad X_{+}\,=\,q^{1/2}EK, \quad X_{-}\,=\,q^{-1/2}FK\}$,
\begin{align}
\cq_{4}\,=\{a+q^{-1}a^*-(1+q^{-1}), c^2, c^{*2}, cc^*, (a-1)c, (a-1)c^*\};
\label{q4}
\\ \nn
\end{align}
\item $\cx_{\cq_{5}}=\{X_{z}\,=\,\frac{\eps_{-}K^{2}-1}{q^{-1}+1}, \quad X_{+}\,=\,q^{1/2}EK, \quad X_{-}\,=\,q^{-1/2}FK\}$,
 \begin{align}
 \cq_{5}\,=\{a-q^{-1}a^*-(1-q^{-1}), c^2, c^{*2}, cc^*, (a-1)c, (a-1)c^*\};
 \label{q5}
 \\ \nn
 \end{align}
\item $\cx_{\cq_{6}}=\{X_{z}\,=\,q(q^2-1)\{FEK^2\,+\,\frac{q^3(K^4-1)}{(q^2-1)^2}\},\quad  X_{+}\,=\,q^{1/2}EK, \quad X_{-}\,=\,q^{-1/2}FK\}$,
\begin{align}
\cq_{6}\,=\{a+q^{-4}a^*-(1+q^{-4}), c^2, c^{*2}, cc^*+(q^3-q)(a-1), (a-1)c, (a-1)c^*\};
\label{q6}
\\ \nn 
\end{align}
\item  $\cx_{\cq_{7}}=\{X_{z}= (1-q^{-2})^{-1}\frac{1-K^4}{1-q^{-2}},\quad X_{+}=q^{1/2}EK, \quad X_{-}=q^{-1/2} FK\},$
\begin{align}
\cq_{7}\,=\{a+q^{-2}a^*-(1+q^{-2}), c^2, c^{*2}, cc^*, (a-1)c, (a-1)c^*\}.
\label{q7}
\\ \nn
\end{align}
\end{enumerate}
An immediate inspection shows that the first order differential calculus $\cq_7$ is the one introduced by Woronowicz \cite{wor87}; calculi (2) and (5) are obtained by the calculi (1) and (4) after mapping $q\to -q$.   The first order calculi associated to  $\cq_1$ and $\cq_{2}$  come as quotients of the four dimensional bicovariant  $4D_{\pm}$ calculi  introduced in \cite{wor89}.

Exterior algebras and differential calculi built over the first order calculi in $\ck$ using the   braidings from \cite{hec2001} present interesting common aspects which can be easily proved by straightforward although long computations that we prefer to omit. We remark that a general proof for these results does not exist, since the braidings we adopted are not canonical. 
By presenting them in the following lines, 
 we also omit to explicitly write any  dependence   (of the bimodules, forms, braidings and antisymmetrisers) 
on index $j=(1,\ldots,7)$  labelling the calculi.  

Given the basis of the quantum tangent space $\cx_{\cq}$, exact one-forms can be written \eqref{ded} as
\beq
\dd x\,=\,\sum_{a}(X_{a}\,\lt \,x)\omega_{a}\qquad\qquad a=\{\pm,z\}
\label{ex1f}
\eeq
with $x\,\in\,\ASU$ on the dual basis of left-invariant one forms. The antilinear hermitian conjugation on $\Gamma_{{\rm L}}$ is 
$\omega_{-}^{*}\,=\,-\omega_{+}, \; \omega_z^{*}\,=\,-\omega_{z}$. The $\ASU$-bimodule is right $\U(1)$-covariant with respect to the natural extension $\delta_{R}^{(1)}:\Gamma\to\Gamma\otimes\ca(\U(1))$
to one forms of the coaction \eqref{cancoa}, set by
\beq
\label{extco}
\delta_{R}^{(1)}(\omega_{z})=\omega_{z}\otimes 1, \qquad\qquad 
\delta_{R}^{(1)}(\omega_{\pm})=\omega_{\pm}\otimes z^{\pm2}.
\eeq
The braiding and its inverse $\sigma, \,\sigma^{-1}:\,\Gamma^{\otimes2}\to\Gamma^{\otimes2}$ are compatible with this $\U(1)$ grading, and have the following spectral decomposition,
\begin{align}
&(1-\sigma)(q^2+\sigma)=0, \nn \\& (1-\sigma^{-1})(q^{-2}+\sigma^{-1})=0
\label{spbra} 
\end{align}
with 
\begin{align}
&\dim\,\ker(1-\sigma^{\pm})\,=6, \nn \\&\dim\,\ker(q^{\pm2}+\sigma^{\pm})\,=3.
\label{degen}
\end{align}
 The antisymmetriser operators  $A_{\sigma^{\pm}}^{(k)}\,:\,\Gamma^{\otimes k}\to\Gamma^{\otimes k}$ they give rise  can be written as
\begin{align}
A^{(2)}_{\sigma^{\pm}}=1-\sigma^{\pm},\qquad\qquad&{\rm on}\,\Gamma^{\otimes 2} \nn \\ 
A^{(3)}_{\sigma^{\pm}}=(1-\sigma^{\pm}_{2})(1-\sigma^{\pm}_1+\sigma^{\pm}_1\sigma^{\pm}_2),\qquad\qquad&{\rm on} \,\Gamma^{\otimes 3} 
\label{anespl}
\end{align}
with $\sigma^{\pm}_1=(\sigma^{\pm}\otimes 1)$ and $\sigma^{\pm}_2=(1\otimes \sigma^{\pm})$, while 
$A^{(k)}_{\sigma^{\pm}}$   are trivial for $k\geq4$. 
They yield  an isomorphic (as left-covariant $\ASU$-bimodules) pair of  exterior algebras ${\rm Range}(A^{(k)}_{\sigma})=\Gamma^{k}_{\sigma}\sim\Gamma^{k}_{\sigma^{-}}={\rm Range}(A^{(k)}_{\sigma^{-}})$ whose dimensions 
coincide to those in the classical setting,  namely $\dim\,\Gamma^{k}_{\sigma^{\pm}}\,=\,3!/(3-k)!$.  A basis of left invariant two forms in $\Gamma_{\sigma}^{2}$ is given by $\{\omega_{-}\wedge\omega_{+}, \omega_{+}\wedge\omega_{z}, \omega_{z}\wedge\omega_{-}\}$; it allows to write the isomorphism above as $\omega_{a}\wedge\omega_{b}\,=\,q^2\omega_{a}\vee\omega_{b}$ (with $a\neq b$) where the symbol $\vee$ clearly represents the wedge product in the exterior algebras  $\Gamma_{\sigma^-}$. 
Given $\vartheta=\omega_{-}\otimes\omega_{+}\otimes\omega_{z}$, left invariant volume forms are then, up to complex numbers,  
\beq
\theta_{\pm}\,=\,A^{(3)}_{\sigma^{\pm}}(\vartheta)
\label{volf}
\eeq
with $\theta_{-}\,=\,q^{-6}\theta_{+}$.

Equipped with the natural graded extension of the exterior differential  in \eqref{ex1f}, these exterior algebras  give isomorphic differential calculi $(\dd, \Gamma_{\sigma})\sim(\dd,\Gamma_{\sigma^{-1}})$. Such differential calculi turn out to be isomorphic to the universal calculus $(\dd, \Gamma_{u})$, since  the relation $\ker\,A^{(k)}_{\sigma^{\pm}}\,=\,\mathcal{S}_{\cq}$ among 2-sided ideals (see section \ref{ss:qs}) holds. 

The action of the antisymmetrisers $A_{\sigma^{\pm}}^{(k)}$ on $\Gamma^{k}_{\sigma^{\pm}}$ is constant. Their  spectral resolution $A^{(k)}_{\sigma^{\pm}}(\phi)\,=\,\lambda_{(k)}^{\pm}\phi$ with $\phi$ a $k$-form yields:
\beq
\lambda_{(2)}^{\pm}\,=\,(1+q^{\pm2}) \qquad\qquad \lambda^{\pm}_{(3)}\,=\,(1+q^{\pm2})(1+q^{\pm2}+q^{\pm4}). 
\label{spans}
\eeq
The isomorphisms $\Gamma^{k}_{\sigma}\sim\Gamma^{k}_{\sigma^{-1}}$ can  then  be written in terms of these spectral resolutions:
\beq
\frac{\omega_{a}\wedge\omega_{b}}{\lambda_{(2)}^{+}}\,=\,
\frac{\omega_{a}\vee\omega_{b}}{\lambda_{(2)}^{-}}, \qquad\qquad
\frac{\theta_+}{\lambda_{(3)}^+}\,=\,\frac{\theta_-}{\lambda_{(3)}^{-}}.
\label{isospec}
\eeq

\begin{rema}
Considering the properties of the so called Drinfeld-Radford-Yetter modules, it is possible to define a braiding $\Psi$ on $\Gamma^{\otimes2}$ as a map which satisifies a braid equation \eqref{br3} on $\Gamma^{\otimes 3}$, using the $\U(1)$-right covariance of the calculi on $\SU$ as explained in \cite{majbo}. We notice that such a braiding has been  considered in \cite{bema, bua} for the case of the Woronowicz' calculus, together with a further braiding $\sigma$ (depending on $\Psi$) obtained in order to construct a meaningful Killing metric, and that they both do not coincide with the braiding we are using in this paper, coming from \cite{hec2001}.

Explicit calculations show, this happens for any of the calculi we consider.  The spectral decompositions  of the braiding  $\Psi$ corresponding   to any of the calculi in $\ck$  presents  $\dim\,\ker(1-\Psi)\,=\,2$
 (see \eqref{degen}). This shows moreover that $\mathcal{S}_{\cq}\nsubseteq\ker(1-\Psi)$: the exterior algebra $\Gamma_{\Psi}$ (built using the braiding $\Psi$) is not a quotient of the universal exterior algebra $\Gamma_{u}$ built over the first order differential calculus $\Gamma$ as described in section \ref{ss:qs}. A direct inspection of section 6 in \cite{bema}
moreover shows the differences between the braiding $\sigma$ on the 3D Woronowicz' calculus used there and the braidings associated to $\ck$ in our approach.
\end{rema}

\section{Hodge operators and symmetric contractions over $\SU$}
\label{hosc:s}

The question is now to exploit how it is possible to suitably translate the classical path described in section \ref{ove:ss} into a quantum path towards the introduction of  a notion of Hodge duality operators and of symmetric contractions on the exterior algebras $\Gamma_{\sigma^{\pm}}$.  We start by introducing an operator which parallels the classical one defined in \eqref{fueq}.

\subsection{Contractions and symmetry}\label{ss:cs}

Since we are interested in Hodge duality operators whose corresponding Laplacians map line bundles elements $\cl_{n}\subset\ASU$ to themselves, we consider the class of non degenerate
$\ASU$-left invariant  and $\U(1)$-right invariant contractions. We define it as the set of maps
$\gm\,:\,\Gamma_{\rm L}\times\Gamma_{\rm L}\,\to\,\IC$,  provided they fullfill the condition $\gm(\omega_{a},\omega_{b})\,=\,0$ if $n_{a}+n_{b}\neq 0$ with $\delta_{R}^{(1)}\,:\,\omega_{j}\mapsto \omega_{j}\otimes z^{n_{j}}$  via \eqref{extco}.
The only non zero coefficients  of the contraction are then (non degeneracy being equivalent to $\alpha\,\beta\,\gamma\,\neq\,0$)
\beq
\label{g1fo}
\gm(\omega_{-},\omega_{+})=\alpha, \qquad\quad \gm(\omega_{+},\omega_{-})=\beta,\qquad\quad \gm(\omega_{z},\omega_{z})=\gamma.
\eeq
This contraction is naturally extended to the left invariant part of $\Gamma_{\sigma^{\pm}}$; recalling the classical expressions \eqref{clac}, \eqref{clacon}, via the actions of the quantum antisymmetrisers $A_{\sigma^{\pm}}^{(k)}$ \eqref{anespl} we set
$$
\gm(\omega_{a_1}\wedge\ldots\wedge\omega_{a_{k}}, \omega_{b_1}\wedge\ldots\wedge\omega_{b_{s}})\,=\,\gm(A^{(k)}_{\sigma}(\omega_{a_1}\otimes\ldots\otimes\omega_{a_{k}}),A^{(s)}_{\sigma}( \omega_{b_1}\otimes\ldots\otimes\omega_{b_{s}}))
$$
together with the obvious analog definition on $\Gamma_{\sigma^{-}}$. 
From  the ordered set of one forms given by $\vartheta\,\in\,\Gamma^{\otimes 3}_{\rm L}$, we define the quantum determinants of the contraction $\gm$ -- with respect to the braidings $\sigma^{\pm}$ -- as
\beq
\label{qdete}
{\rm det}_{\sigma^{\pm}}\gm\,=\,\frac{1}{\lambda_{(3)}^{\pm}}\,\gm(\theta_{\pm},\theta_{\pm}),
\eeq
reading
${\rm det}_{\sigma}\gm\,=\,q^{6}{\rm det}_{\sigma^-}\gm$. 
We set then the hermitian  volume forms  $\mu_{\pm}\,=\,m_{\pm}\theta_{\pm}\,=\,\mu_{\pm}^*$ from  \eqref{volf} with $m_{\pm} \,\in\,\IR$, and generalising the classical  \eqref{fueq} to the quantum setting,
  the linear $\ASU$- linear operators $S_{\sigma^{\pm}}:\Gamma^{k}_{\sigma^{\pm}}\to\Gamma^{3-k}_{\sigma^{\pm}}$ as  
\beq
\label{defS}
S_{\sigma^{\pm}}(\omega)\,=\,\frac{1}{\lambda^{\mp}_{(k)}}\,\gm(\omega, \mu_{\pm})
\eeq
on a left-invariant basis; the modulus of the scale factors of the volume are chosen by  $S^2_{\sigma^{\pm}}(1)={\rm sgn(\det}_{\sigma^{\pm}}\gm)$.   Corresponding to these operators we introduce  sesquilinear $\ASU$-left invariant scalar products by
\begin{align}
&\{\omega, \omega^{\prime}\}_{\sigma}\,=\,\int_{\mu_{+}}\omega^*\wedge S_{\sigma}(\omega^{\prime}), \nn \\
&\{\omega, \omega^{\prime}\}_{\sigma^{-}}\,=\,\int_{\mu_{-}}\omega^*\vee S_{\sigma^{-}}(\omega^{\prime}).
\label{sdpm}
\end{align} 
The integral on $\Gamma^3_{\sigma^{\pm}}$ is defined in terms of the Haar functional $h$ by $\int_{\mu_{\pm}}x\,\mu_{\pm}\,=\,h(x)$ for $x\in\,\ASU$.  
The isomorphisms \eqref{isospec} allow to prove the following relations:
\begin{align}
S_{\sigma^-}(1)&=\,\left(\frac{m_-}{m_+}\right)\left(\frac{\lambda^-_{(3)}}{\lambda^+_{(3)}}\right) S_{\sigma}(1),
\nn \\
S_{\sigma^-}(\omega)&=\,\left(\frac{m_-}{m_+}\right)\left(\frac{\lambda^-_{(3)}}{\lambda^+_{(3)}}\right) S_{\sigma}(\omega),
\nn \\ 
S_{\sigma^-}(\phi)&=\,\left(\frac{m_-}{m_+}\right)\left(\frac{\lambda^-_{(2)}}{\lambda^+_{(2)}}\right)\left(\frac{\lambda^-_{(3)}}{\lambda^+_{(3)}}\right) S_{\sigma}(\phi), \nn \\
S_{\sigma^-}(\theta)&=\,\left(\frac{m_-}{m_+}\right)\left(\frac{\lambda^-_{(3)}}{\lambda^+_{(3)}}\right)^2 S_{\sigma}(\theta),
\label{ssm}
\end{align}
for any 1-form $\omega$, 2-form $\phi$, 3-form $\theta$. For the common scale factor one has
\beq
\left(\frac{m_{-}}{m_+}\right)^2\,=\,\left(\frac{\lambda^{+}_{(3)}}{\lambda^{-}_{(3)}}\right)^3
\label{sfso}
\eeq
while, for the scalar products,
\begin{align}
\{\omega,\omega^{\prime}\}_{\sigma^-}&=\,
\left(\frac{\lambda^+_{(2)}}{\lambda^-_{(2)}}\right)\left(\frac{\lambda^-_{(3)}}{\lambda^+_{(3)}}\right) 
\{\omega,\omega^{\prime}\}_{\sigma}, \nn \\
\{\phi,\phi^{\prime}\}_{\sigma^-}&=\,
\left(\frac{\lambda^-_{(3)}}{\lambda^+_{(3)}}\right) 
\{\phi,\phi^{\prime}\}_{\sigma}, \nn \\
\{\theta,\theta^{\prime}\}_{\sigma^-}&=\,
\left(\frac{\lambda^-_{(3)}}{\lambda^+_{(3)}}\right) 
\{\theta,\theta^{\prime}\}_{\sigma}
\label{pscb}
 \end{align}
for any pair $\omega, \omega^{\prime}$ of 1-forms, $\phi, \phi^{\prime}$ of 2-forms, $\theta,\theta^{\prime}$ of 3-forms.



\subsection{Scalar products and duality operators}
\label{ss:fdh}

Before analysing the spectral properties of the operators above, we use 
 the spectral resolution of the antisymmetrisers further and introduce 
via the contraction map  another sesquilinear $\ASU$-left invariant scalar product on the exterior algebras $\Gamma_{\sigma^{\pm}}$ by
\beq
\hs{x\,\omega}{x^{\prime}\omega^{\prime}}_{\sigma^{\pm}}\,=\,h(x^{*}x^{\prime})\,\frac{1}{\lambda_{(k)}^{\pm}}\,\gm(\omega^*,\omega^{\prime}),
\label{scapro}
\eeq
where $\omega, \omega^{\prime}$ are left-invariant forms and $h(x^*x^{\prime})$ is again the action of the Haar functional $h$ with $x,x^{\prime}\,\in\,\ASU$.  Recalling the classical definition \eqref{claTop},  
it is now natural to set the operators\footnote{This operator is the one described in \cite{funingeo} for the Woronowicz calculus.} $T_{\sigma^{\pm}}:\,\Gamma^{ k}_{\sigma^{\pm}}\,\to\,\Gamma^{ 3-k}_{\sigma^{\pm}}$ by:
\begin{align}
&\hs{x\,\omega}{x^{\prime}\omega^{\prime}}_{\sigma}\,=\,\int_{\mu_{+}}(x\,\omega)^*\wedge\,T_{\sigma}(x^{\prime}\omega^{\prime}),\nn \\
&\hs{x\,\omega}{x^{\prime}\omega^{\prime}}_{\sigma^{-}}\,=\,\int_{\mu_{-}}(x\,\omega)^*\vee\,T_{\sigma^{-}}(x^{\prime}\omega^{\prime}).
\label{dTo}
\end{align}
Provided the contraction $\gm$ is non degenerate, such operators $T_{\sigma^{\pm}}$ are well-defined, bijective and left $\ASU$-linear \cite{kmt}. 
One immediately has
\begin{align}
T_{\sigma^{\pm}}(1)&=\mu_{\pm}, \nn \\
T_{\sigma^{\pm}}(\mu_{\pm})&=\hs{\mu_{\pm}}{\mu_{\pm}}_{\sigma^{\pm}}\,=\,m^2_{\pm}\,{\rm det}_{\sigma^{\pm}}\gm.
\label{boh}
\end{align}
From \eqref{g1fo}, the scalar product \eqref{scapro} reads on left-invariant 1-forms (we omit the subscripts $\sigma^{\pm}$ since they coincide)
\beq
\label{se1fo}
\hs{\omega_{-}}{\omega_{-}}=-\beta, \qquad\qquad \hs{\omega_{+}}{\omega_{+}}=-\alpha,\qquad\qquad \hs{\omega_{z}}{\omega_{z}}=-\gamma.
\eeq
If we assume the natural normalization condition $T^2_{\sigma^{\pm}}(1)\,=\,{\rm sgn}({\rm det}_{\sigma^{\pm}}\gm)$, it is easy to prove that 
\begin{align}
\hs{\phi}{\phi^{\prime}}_{\sigma^-}&=\,\left(\frac{\lambda_{(2)}^-}{\lambda_{(2)}^{+}}\right)^3\, 
\hs{\phi}{\phi^{\prime}}_{\sigma}, \nn \\
\hs{\theta}{\theta^{\prime}}_{\sigma^{-}}\,&=\,\left(\frac{\lambda^-_{(3)}}{\lambda^{+}_{(3)}}\right)^3\,\hs{\theta}{\theta^{\prime}}_{\sigma}
\label{eqmp}
\end{align}
which are the counterparts of the \eqref{pscb} (with which they share the same notations) for the operators $T_{\sigma^{\pm}}$ while,  
from \eqref{eqmp},
we have
\begin{align}
&\lambda_{(3)}^-m_-^2\,=
\lambda_{(3)}^+m_{+}^2
\nn \\
&T_{\sigma^{-}}(1)\,=\,\left(\frac{m_{-}}{m_{+}}\right)\left(\frac{\lambda_{(3)}^-}{\lambda_{(3)}^+}\right)T_{\sigma}(1)
,
\nn \\
&T_{\sigma^-}(\omega)\,=\,\left(\frac{m_-}{m_+}\right)\left(\frac{\lambda_{(2)}^-}{\lambda_{(2)}^+}\right)T_{\sigma}(\omega)
,
\nn \\
&T_{\sigma^-}(\phi)\,=\,\left(\frac{m_-}{m_+}\right)T_{\sigma}(\phi)
,
\nn \\
&T_{\sigma^{-}}(\theta)\,=\,\left(\frac{m_-}{m_+}\right)T_{\sigma}(\theta),
\label{repro}
\end{align}
which are the counterparts of the previous \eqref{ssm}, \eqref{sfso}.
It is now evident that the operators $(T_{\sigma^{\pm}},\,S_{\sigma^{\pm}})$ differ, since the  sesquilinear products \eqref{scapro}, \eqref{sdpm}  differ  and that they coincide only in the classical limit; one can easily for example check 
 on 3-forms that  
\beq
\{\theta_{\pm},\theta_{\pm}\}_{\sigma^{\pm}}\,=\,\frac{1}{\lambda_{(3)}^{\mp}}\,\gm(\theta_{\pm},\theta_{\pm})\,=\,\frac{\lambda^{\pm}_{(3)}}{\lambda^{\mp}_{(3)}}\,\hs{\theta_{\pm}}{\theta_{\pm}}_{\sigma^{\pm}}.
\label{covo}
\eeq
In order to better understand the differences between the operators $S_{\sigma^{\pm}}, \,T_{\sigma^{\pm}}$ as well as their similarities, we explicitly present a deeper analysis in the case of the Woronowicz calculus $(\dd, \Gamma_{\sigma^{\pm}}^{\mathcal{W}})$, that we use as an example. 

\subsection{The guiding example}
On the Woronowicz first order differential calculus  the braiding reads  
\begin{align}
\sigma(\omega_{a}\otimes\omega_{a})=\omega_{a}\otimes\omega_{a},&\qquad\qquad a=\pm,z
\nn \\
\sigma(\omega_{-}\otimes\omega_{+})=(1-q^{2})\omega_{-}\otimes\omega_{+}+q^{-2}\omega_{+}\otimes\omega_{-}, & \qquad\qquad \sigma(\omega_{+}\otimes\omega_{-})=q^4\omega_{-}\otimes\omega_{+}, \nn \\
\sigma(\omega_{-}\otimes\omega_{z})=(1-q^{2})\omega_{-}\otimes\omega_{z}+q^{-4}\omega_{z}\otimes\omega_{-}, &
\qquad\qquad \sigma(\omega_{z}\otimes\omega_{-})=q^6\omega_{-}\otimes\omega_{z}, \nn \\
\sigma(\omega_{z}\otimes\omega_{+})=(1-q^{2})\omega_{z}\otimes\omega_{+}+q^{-4}\omega_{+}\otimes\omega_{z},& \qquad\qquad
\sigma(\omega_{+}\otimes\omega_{z})=q^6\omega_{z}\otimes\omega_{+},
\label{braiwor}
\end{align}
so that the wedge product  on the exterior algebra satisfies
$\omega_{a}\wedge\omega_{a}=0, \quad (a=\pm,z)$ with  
\beq
\label{commc3wo}
\omega_{-}\wedge\omega_{+}+q^{-2}\omega_{+}\wedge\omega_{-}=0,
\qquad\qquad
\omega_{z}\wedge\omega_{\mp}+q^{\pm4}\omega_{-}\wedge\omega_{z}=0.
\eeq 
From \eqref{volf} and \eqref{isospec} it is
\begin{align}
&\theta\,=\,q^4\,\om\otimes(\op\otimes\oz\,-\,q^6\,\oz\otimes\op)\,+\,q^{-6}\,\op\otimes(\oz\otimes\om\,-\,q^6\om\otimes\oz)
\nonumber \\ &\qquad\qquad\qquad\qquad\qquad\qquad\qquad\qquad\qquad\qquad \,+\,q^{4}\,\oz\otimes(\om\otimes\op\,-\,q^{-4}\op\otimes\om)
 \label{volw}
 \end{align}
 and 
\beq
\label{gmwo}
\gm(\theta_{+},\theta_{+})\,=\,-6\,q^4\,\alpha\,\beta\,\gamma.
\eeq
Together with \eqref{boh}, one has 
\begin{align}
T_{\sigma}(\omega_{-})\,=\,q^{-2}m_+\,\hs{\om}{\om}\,\om\wedge\oz, &\qquad\qquad T_{\sigma}(\om\wedge\oz)\,=\,q^{-2}m_{+}\hs{\om\wedge\oz}{\om\wedge\oz}_{\sigma}\,\om \nn \\
T_{\sigma}(\omega_{+})\,=\,-m_+\,\hs{\op}{\op}\,\op\wedge\oz, &\qquad\qquad T_{\sigma}(\op\wedge\oz)\,=\,-m_{+}\hs{\op\wedge\oz}{\op\wedge\oz}_{\sigma}\,\op\, \nn \\
T_{\sigma}(\omega_{z})\,=\,-m_+\,\hs{\oz}{\oz}\,\om\wedge\op,&\qquad\qquad T_{\sigma}(\om\wedge\op)\,=\,-m_{+}\hs{\om\wedge\op}{\om\wedge\op}_{\sigma}\,\oz\,
\label{Texpw1}
\end{align}
while the relations \eqref{repro} give the action of $T_{\sigma^-}$.  
The expressions above depend only on the wedge products relations \eqref{commc3wo}; they are valid for any choice of a non degenerate scalar product on the space of left-invariant  1-, 2- and 3-forms.  Since from the expression \eqref{sdpm} we see that the scalar product $\{~,~\}_{\sigma^{\pm}}$ characterises the operators $S_{\sigma^{\pm}}$ in the same way the scalar product $\hs{~}{~}_{\sigma^{\pm}}$ characterises the operators $T_{\sigma^{\pm}}$ (see \eqref{dTo}), it is immediate to recover that the action of the operators $S_{\sigma^{\pm}}$ can be written from the  action of the operators $T_{\sigma^{\pm}}$ after the mapping  $\hs{~}{~}_{\sigma^{\pm}}\,\mapsto\,\{~,~\}_{\sigma^{\pm}}$ on each space $\Gamma^{k}_{\sigma^{\pm}}$.

But we want to show that a deeper analogy exists.  It is clear from the structure of the braiding that the scalar product \eqref{scapro} on left-invariant $k$-forms ($k\,=\,2,3$) is a $k$-order polynomial in the first order terms \eqref{se1fo}. This means that the definition \eqref{scapro} amounts to set a specific choice for the extension to higher order forms of the scalar product $\hs{~}{~}_{\sigma}$ on 1-forms. For the example we are considering, we calculate 
\begin{align}
\hs{\om\wedge\op}{\om\wedge\op}_{\sigma}\,=\,2\,\hs{\om}{\om}\,\hs{\op}{\op}/\lambda^+_{(2)}, \nn &\\
\hs{\om\wedge\oz}{\om\wedge\oz}_{\sigma}\,=\,2\,q^{-2}\hs{\om}{\om}\,\hs{\oz}{\oz}/\lambda^+_{(2)}, &\nn \\
\hs{\op\wedge\oz}{\op\wedge\oz}_{\sigma}\,=\,2\,q^{6}\hs{\op}{\op}\,\hs{\oz}{\oz}/\lambda^+_{(2)}, &
\label{exwo2}
\end{align}
on 2-forms, and (see \eqref{gmwo})
\beq
\hs{\theta_{+}}{\theta_{+}}_{\sigma}\,=\,\frac{6\,q^4}{\lambda_{(3)}^{+}}\,\hs{\om}{\om}\,\hs{\op}{\op}\,\hs{\oz}{\oz}.
\label{equ3wo}
\eeq
on the volume form. Concerning the  scalar product \eqref{sdpm},   we start by computing that we have
\begin{align}
&\{\om,\om\}_{\sigma}\,=\,-\alpha,\qquad\qquad
\{\op,\op\}_{\sigma}\,=\,-q^4\beta,\qquad\qquad
\{\oz,\oz\}_{\sigma}\,=\,-q^2\gamma \nn \\
&\{\om,\om\}_{\sigma^-}\,=\,-q^{-4}\alpha,\qquad\qquad
\{\op,\op\}_{\sigma^-}\,=\,-\beta,\qquad\qquad
\{\oz,\oz\}_{\sigma^-}\,=\,-q^{-2}\gamma 
\label{se1}
\end{align}
on 1-forms. 
 The next step is to understand how this scalar product on higher order forms can be written in terms of the scalar products among 1-forms. It turns then out that we can write:
\begin{align}
&\{\om\wedge\op,\om\wedge\op\}_{\sigma}\,=\,2\,\{\om,\om\}_{\sigma}\,\{\op,\op\}_{\sigma}/\lambda^{+}_{(2)}, \nn \\&
\{\om\wedge\oz,\om\wedge\oz\}_{\sigma}\,=\,2\,q^{-2}\{\om,\om\}_{\sigma}\,\{\oz,\oz\}_{\sigma}/\lambda^+_{(2)}, \nn \\&
\{\op\wedge\oz,\op\wedge\oz\}_{\sigma}\,=\,2\,q^{6}\{\op,\op\}_{\sigma}\,\{\oz,\oz\}_{\sigma}/\lambda^+_{(2)}, 
\label{exwo2c}
\end{align}
 and (see \eqref{gmwo})
\beq
\{\theta_{+},\theta_{+}\}_{\sigma}\,=\,\frac{6\,q^4}{\lambda_{(3)}^{+}}\,\{\om,\om\}_{\sigma}\,\{\op,\op\}_{\sigma}\,\{\oz,\oz\}_{\sigma}.
\label{equ3woc}
\eeq
on the volume form. The comparison is immediate: with respect to first order coefficients, the equations \eqref{exwo2c} have the same structure of  the equations \eqref{exwo2}, the  equation \eqref{equ3woc} analogously has the same structure of the equation \eqref{equ3wo}. This equivalence holds also if we consider  the scalar product defined by  the operator $S_{\sigma^{-}}$, as one may easily infers from 
\begin{align}
\frac{\{\omega_a\vee\omega_b,\omega_a\vee\omega_b\}_{\sigma^{-}}}{\{\omega_a,\omega_a\}_{\sigma^-}\,\{\omega_b,\omega_b\}_{\sigma^{-}}}
\,=\,\frac{\lambda_{(2)}^-}{\lambda_{(2)}^+}\,\frac{\{\omega_a\wedge\omega_b,\omega_a\wedge\omega_b\}_{\sigma}}{\{\omega_a,\omega_a\}_{\sigma}\,\{\omega_b,\omega_b\}_{\sigma}}
\nn \\
\frac{\{\theta_-,\theta_-\}_{\sigma^{-}}}{\{\omega_-,\omega_-\}_{\sigma^-}\,\{\omega_+,\omega_+\}_{\sigma^{-}}\,\{\oz,\oz\}_{\sigma^{-}}}
\,=\,\frac{\lambda_{(3)}^-}{\lambda_{(3)}^+}\,\frac{\{\theta_+,\theta_+\}_{\sigma}}{\{\omega_-,\omega_-\}_{\sigma}\,\{\omega_+,\omega_+\}_{\sigma}\,\{\oz,\oz\}_{\sigma}},
\label{inte}
\end{align}
which parallel, for the specific example of the Woronowicz calculus, the relations \eqref{eqmp}. 

The collection of these results prove that the scalar products $\{~,~\}_{\sigma^{\pm}}$ among higher order (left-invariant) forms   can be obtained from $\hs{~}{~}_{\sigma^{\pm}}$ via the replacements of  the 1-forms terms, as (restoring only for this expression the index $\sigma^{\pm}$ to the $\hs{~}{~}$ scalar product between 1-forms)
\begin{align}
\hs{\omega_a}{\omega_b}_{\sigma}\quad&\leftrightarrow\quad\{\omega_a, \omega_b\}_{\sigma}\,=\,\frac{1}{\lambda^{-}_{(k)}}\,\int_{\mu_{+}}\omega_{a}^*\wedge \gm(\omega_b, \mu_{+}), \nn \\
\hs{\omega_a}{\omega_b}_{\sigma^{-1}}\quad&\leftrightarrow\quad\{\omega_a, \omega_b\}_{\sigma^{-}}\,=\,\frac{1}{\lambda^{+}_{(k)}}\,\int_{\mu_{-}}\omega_a^*\vee \gm(\omega_b, \mu_{-});
\label{sde1}
\end{align}  
the comparison between \eqref{sdpm} and \eqref{dTo} convinces that the arrows in \eqref{sde1} give the action of $T_{\sigma^{\pm}}$ from the action of $S_{\sigma^{\pm}}$ and viceversa. 


\subsection{A shared pattern}
We come now to the crucial point of our  analysis. Does the  equivalence -- described above for our guiding example --  between the scalar products defined by \eqref{scapro} and \eqref{sdpm} hold also for the other calculi in $\ck$  on $\SU$? Do  -- for any fixed calculus in $\ck$ --   these two families share the same pattern, once we look at the scalar products among higher order forms  as polynomials over  their first order coefficients? The answer is positive, and can be proved by straightforward but long explicit computations, since we miss a general theory for the  set of  non canonical braidings \cite{hec2001} we consider.  Relations \eqref{inte} hold for any of the calculi in $\ck$;  the mappings in \eqref{sde1} allow to obtain the action of the operators $T_{\sigma^{\pm}}$ from that of the operators $S_{\sigma^{\pm}}$ and viceversa.

\subsection{Symmetric and real contractions}
\label{shs}

We now use these operators to introduce a notion of symmetry and reality for the contraction $\gm$. Let us define the contraction $\gm$:
\begin{enumerate}[(i)]
\item {\it $S_{\sigma^{\pm}}$-symmetric} (resp. {\it $T_{\sigma^{\pm}}$-symmetric)}, provided the operators $S^2_{\sigma^{\pm}}$ (resp. $T_{\sigma^{\pm}}^2$)  have the same degeneracy of the antisymmetrisers, namely is their action on one forms constant;
\item  {\it $S_{\sigma^{\pm}}$-real} (resp. {\it $T_{\sigma^{\pm}}$-real}), provided the relations $S_{\sigma^{\pm}}(\omega_{a}^*)\,=\,(S_{\sigma^{\pm}}(\omega_{a}))^*$ (resp. $T_{\sigma^{\pm}}(\omega_{a}^*)\,=\,(T_{\sigma^{\pm}}(\omega_{a}))^*$)  on any left-invariant one form hold.
\end{enumerate}
Denote by $G_S^{\pm}$ (resp. $G_{T}^{\pm}$)   the set of  real and symmetric contractions: we consider then the corresponding dualities $S_{\sigma^{\pm}}$ (resp. $T_{\sigma^{\pm}}$)  as Hodge operators. 
The relations  \eqref{repro} enable to prove that $G_{S}^+=G_{S}^-$, while  the relations  \eqref{inte}  give also $G_T^+=G_{T}^-$ (from now on we shall then denote these sets by $G_S,\,G_T$).  The requirements of reality and symmetry clearly amount to constraint the parameters $\alpha, \,\beta, \gamma$; such sets are not void, and do not  coincide, i.e. $G_{T}\,\neq\, G_S$. 

It is easy to compute, for the Woronowicz calculus,
\begin{align}
&\gm\,\in\,G_{T}\qquad\Leftrightarrow\qquad\{\beta\,=\,q^6\alpha\,\in\,\IR,\;\gamma\,\in\,\IR\} \nn \\
&\gm\,\in\,G_{S}\qquad\Leftrightarrow\qquad\{\beta\,=\,q^{10}\alpha\,\in\,\IR,\;\gamma\,\in\,\IR\} \label{wosym}
 \end{align}
Following the long analysis of the two scalar products and their corresponding Hodge duality operators, comes naturally  that  these two sets of constraints are equivalent, if they are written in terms of the corresponding scalar products on 1-forms, namely  
 \begin{align}
&\gm\,\in\,G_{T}\qquad\Leftrightarrow\qquad\{\hs{\om}{\om}\,=\,q^6\hs{\op}{\op}\,\in\,\IR,\;\hs{\oz}{\oz}\,\in\,\IR\} \nn \\
&\gm\,\in\,G_{S}\qquad\Leftrightarrow\qquad\{\{\om,\om\}_{\sigma^{\pm}}\,=\,q^{6}\{\op,\op\}_{\sigma^{\pm}}\,\in\,\IR,\;\{\oz,\oz\}_{\sigma^{\pm}}\,\in\,\IR\}. \label{wosym1}
 \end{align}
 Such an equivalence is again  given by \eqref{sde1}.  We close this part by saying that the equivalences given in \eqref{wosym1} and \eqref{laeq} hold for any of the other six calculi \eqref{q1}-\eqref{q6} in $\ck$.  This means, that what we have actually  introduced is (for any of the calculi in $\ck$ on $\SU$) a notion of real and symmetric tensor  over the  vector space $\mathcal{X}^{\otimes 2}\,=\,\mathcal{X}\otimes_{\IC}\mathcal{X}$, 
which we write as
\beq
\gm\,=\,\sum_{a,b\,=\,\pm,z}\gm_{ab}\,X_{a}\otimes X_{b} 
 \label{metr}
 \eeq
with $\gm_{ab}\,=\,\gm(\omega_{a},\omega_{b})$ belonging either to 
$G_{S}$ or $G_{T}$.  We notice that no compelling reason  at this level of mathematical analysis allows to select the notion of $T_{\sigma^{\pm}}$-symmetry with respect to the one of $S_{\sigma^{\pm}}$-symmetry, nor to select one of the calculi we have considered.

 What differs (only for $q\neq1$) are the spectra of the Laplacians associated to the Hodge operators $T_{\sigma^{\pm}},\,S_{\sigma^{\pm}}$ and the symmetric tensor \eqref{metr}, which turn out to be
\begin{align}
\Box^{(T)}x\,&=\,{\rm sgn(det}_{\sigma^{\pm}} \gm)\,\sum_{a,b\,=\,\pm,z}\left[\hs{\omega_a^*}{\omega_{b}}\,X_{a}X_{b}\,\lt\,x\right] \nn \\ 
\Box^{(S)}x\,&=\,{\rm sgn(det}_{\sigma^{\pm}} \gm)\,\sum_{a,b\,=\,\pm,z}\left[\{\omega_a^*,\omega_{b}\}_{\sigma^{\pm}}\,X_{a}X_{b}\,\lt\,x\right]. \label{laeq} 
\end{align}
It is interesting to notice that such Laplacians depend on the scalar products associated to the symmetric and real tensor \eqref{metr}.

In the next part we present, for each of these calculi, the sufficient ingredients to build the isomorphic exterior algebras and the Hodge operator $S_{\sigma}$ together with its corresponding class of real and symmetric contractions (those which are $\ASU$-left invariant and $\U(1)$-right  coinvariant). As we saw, this is enough to construct the operators $S_{\sigma^-}, \,T_{\sigma^{\pm}}$.

\subsection{Hodge operators}
\label{ss:Ho}
We follow the numbering in section \ref{ss:sesu}; we assume the non degeneracy of the contraction $\gm$ in \eqref{g1fo}, that is $\alpha\,\beta\,\gamma\,\neq\,0$.

\bigskip

\begin{enumerate}[(1)]

\item
Given the quantum tangent space $\mathcal{X}_{\cq_1}$ the exact one forms are
\begin{align*}
\dd a=-q\,c^*\omega_{+}\,+\,a\,\omega_{z},& \qquad\qquad \dd c=a^*\omega_{+}\,+\,c\,\omega_{z}, 
\\
\dd a^*=c\,\omega_-\,-\,q^{-1}a^*\omega_{z},& \qquad\qquad \dd c^*=-q^{-1}a\,\omega_{-}\,-\,q^{-1}c^*\omega_z;
\end{align*}
and the braiding reads
\begin{align*}
\sigma(\omega_{\pm}\otimes\omega_{\pm})=\omega_{\pm}\otimes\omega_{\pm},& \\
\sigma(\omega_z\otimes\omega_z)=\omega_z\otimes\omega_z\,+\,\frac{q(1-q)}{1+q^{-1}}(\omega_{+}\otimes\omega_-\,-\,\omega_{-}\otimes\omega_{+}),&
\\
\sigma(\omega_-\otimes\omega_+)=q^2\omega_+\otimes\omega_-\,+\,(1-q^2)\omega_-\otimes\omega_+,&\qquad\qquad\sigma(\omega_{+}\otimes\omega_-)=\omega_-\otimes\omega_+ \\
\sigma(\om\ot\oz)=q^2\oz\ot\om\,+\,(1-q^2)\om\ot\oz,& \qquad\qquad \sigma(\oz\ot\om)=\om\ot\oz
\\
\sigma(\oz\ot\op)=q^2\op\ot\oz\,+\,(1-q^2)\oz\ot\op,& \qquad\qquad \sigma(\op\ot\oz)=\oz\ot\op.
\end{align*}
The hermitian structure over left-invariant two forms and the wedge product antisymmetry are
 \begin{align*}
 (\omega_-\wedge\omega_+)^*&=-\omega_-\wedge\omega_{+}=q^{2}\omega_{+}\wedge\omega_{-}, \\
 (\omega_-\wedge\omega_z)^*&=-\omega_z\wedge\omega_{+}=q^{2}\omega_{+}\wedge\omega_{z}, \\
(\omega_+\wedge\omega_z)^*&=-\omega_z\wedge\omega_{-}=q^{-2}\omega_{-}\wedge\omega_{z};
\end{align*}
the volume form $\theta_+$ turns out to be a multiple of the classical one, namely the one we would obtain if  the braiding were the  classical flip,
\beq
\label{vol1}
\theta_+=q^4(\omega_-\otimes(\omega_{+}\otimes\omega_{z}\,-
\omega_z\otimes\omega_{+})\,+\,\op\otimes(\omega_z\otimes\omega_{-}\,-\,\om\ot\oz) 
\,+\,\oz\ot(\om\ot\op\,-\,\op\ot\om)),
\eeq
while for the normalisations of the Hodge operators one needs  
$\gm(\theta_+,\theta_+)=\,-\,6q^8(\alpha\beta\gamma)$. It is
$$
\{\om,\om\}_{\sigma}\,=\,-\alpha, \qquad
\{\op,\op\}_{\sigma}\,=\,-q^4\,\beta, \qquad
\{\oz,\oz\}_{\sigma}\,=\,-q^{2}\,\gamma.
$$
A contraction \eqref{g1fo} is $S_{\sigma^{\pm}}$-real and symmetric 
\beq 
\gm\,\in\,G_{S}\qquad\Leftrightarrow\qquad\{\alpha\,=\,\beta\,\in\,\IR,\, \gamma\,\in\,\IR\},
\label{coq1}
\eeq 
while $\gm\,\in\,G_{T}\Leftrightarrow\{\alpha\,=\,q^4\,\beta\,\in\,\IR,\, \gamma\,\in\,\IR\}$.
As corresponding Hodge operator we have:
\begin{align}
S_{\sigma}(1)=\mu_+, &\qquad\qquad S_{\sigma}(\mu_+)=-{\rm sgn}(\gamma), \nn \\
S_{\sigma}(\om)=m_+\,q^2\{\om,\om\}_{\sigma}\,\om\wedge\oz, &\qquad\qquad S_{\sigma}(\om\wedge\oz)=2m(q^6/\lambda_{(2)}^+)\{\om,\om\}_{\sigma}\,\{\oz,\oz\}_{\sigma}\,\om \nn \\
S_{\sigma}(\op)=-m_{+}\,\{\op,\op\}_{\sigma}\,\op\wedge\oz, &\qquad\qquad S_{\sigma}(\op\wedge\oz)=-2m_+(1/\lambda_{(2)}^+)\{\op,\op\}_{\sigma}\,\{\oz,\oz\}_{\sigma}\,\op \nn \\
S_{\sigma}(\oz)=-m_{+}\,\{\oz,\oz\}_{\sigma}\,\om\wedge\op, &\qquad\qquad S_{\sigma}(\om\wedge\op)=-2m(q^4/\lambda_{(2)}^+)\{\om,\om\}_{\sigma}\,\{\op,\op\}_{\sigma}\,\oz 
\label{Hoq1}
\end{align}
with a normalisation condition $m_+^2\,{\rm det}_{\sigma}\gm\,=\,-{\rm sgn} \gamma$. 
\bigskip 
\bigskip
\item
As we  have already noticed, the structure of the exterior algebra corresponding to this calculus is obtained by that corresponding to the  previous one by mapping  $q\,\to\,-q$. Since the relations \eqref{Hoq1} (and then the \eqref{coq1}) are invariant by this reflection, the Hodge duality we obtain is exactly the previous one.  

\bigskip 
\bigskip
\item
Given the quantum tangent space $\mathcal{X}_{\cq_3}$ the exact one forms are
\begin{align*}
\dd a=-q\,c^*\omega_{+}\,+\,a\,\omega_{z},& \qquad\qquad \dd c=a^*\omega_{+}\,+\,c\,\omega_{z}, 
\\
\dd a^*=c\,\omega_-\,-\,q^{-2}a^*\omega_{z},& \qquad\qquad \dd c^*=-q^{-1}a\,\omega_{-}\,-\,q^{-2}c^*\omega_z;
\end{align*}
with a braiding:
\begin{align*}
\sigma(\omega_{\pm}\otimes\omega_{\pm})=\omega_{\pm}\otimes\omega_{\pm},& \\
\sigma(\omega_z\otimes\omega_z)=\omega_z\otimes\omega_z\,+\,(q^2-1)(\omega_{-}\otimes\omega_+\,-\,q^4\omega_{+}\otimes\omega_{-}),&
\\
\sigma(\omega_-\otimes\omega_+)=q^6\omega_+\otimes\omega_-\,+\,(1-q^2)\omega_-\otimes\omega_+,&\qquad\qquad\sigma(\omega_{+}\otimes\omega_-)=q^{-4}\omega_-\otimes\omega_+ \\
\sigma(\om\ot\oz)=q^4\oz\ot\om\,+\,(1-q^2)\om\ot\oz,& \qquad\qquad \sigma(\oz\ot\om)=q^{-2}\om\ot\oz
\\
\sigma(\oz\ot\op)=q^4\op\ot\oz\,+\,(1-q^2)\oz\ot\op,& \qquad\qquad \sigma(\op\ot\oz)=q^{-2}\oz\ot\op.
\end{align*}
The hermitian structure over left-invariant two forms is
 \begin{align*}
 (\omega_-\wedge\omega_+)^*&=-\om\wedge\op=q^{6}\omega_{+}\wedge\omega_{-}, \\
 (\omega_-\wedge\omega_z)^*&=-\oz\wedge\op=q^{4}\omega_{+}\wedge\omega_{z}, \\
(\omega_+\wedge\omega_z)^*&=-\oz\wedge\om=q^{-4}\omega_{-}\wedge\omega_{z};
\end{align*}
the volume form $\theta_{+}$  is 
\beq
\label{vol3}
\theta_+=q^4\omega_-\otimes(\omega_{+}\otimes\omega_{z}\,-\,q^{-2}
\omega_z\otimes\omega_{+})\,-\,q^8\op\otimes(\omega_-\otimes\omega_{z}\,-\,q^2\oz\ot\om) 
\,+\,q^4\oz\ot(\om\ot\op\,-\,q^4\op\ot\om),
\eeq
so to have $\gm(\theta_{+},\theta_{+})\,=\,-\,6q^{12}(\alpha\beta\gamma).$ 
It is
$$
\{\om,\om\}_{\sigma}\,=\,-\alpha, \qquad
\{\op,\op\}_{\sigma}\,=\,-q^4\,\beta, \qquad
\{\oz,\oz\}_{\sigma}\,=\,-q^{2}\,\gamma.
$$
A contraction \eqref{g1fo} is $S_{\sigma^{\pm}}$-real and symmetric 
\beq 
\gm\,\in\,G_{S}\qquad\Leftrightarrow\qquad\{q^6\,\alpha\,=\,\beta\,\in\,\IR,\, \gamma\,\in\,\IR\};
\label{coq3}
\eeq 
while the conditions of $T_{\sigma}$-reality and symmetry for the same contraction are
$\gm\,\in\,G_{T}\Leftrightarrow\{\alpha\,=\,q^{10}\beta\,\in\,\IR,\, \gamma\,\in\,\IR\}$.
The  Hodge operator is:
\begin{align}
S_{\sigma}(1)=\mu_+, &\qquad\qquad S_{\sigma}(\mu_+)=-{\rm sgn}(\gamma), \nn \\
S_{\sigma}(\om)=m_+\,q^6\{\om,\om\}_{\sigma}\,\om\wedge\oz, &\qquad\qquad S_{\sigma}(\om\wedge\oz)=2m_+(q^{12}/\lambda_{(2)}^+)\{\om,\om\}_{\sigma}\,\{\oz,\oz\}_{\sigma}\,\om \nn \\
S_{\sigma}(\op)=-m_+\,\{\op,\op\}_{\sigma}\,\op\wedge\oz, &\qquad\qquad S_{\sigma}(\op\wedge\oz)=-2m_{+}(q^{-2}/\lambda_{(2)}^+)\{\op,\op\}_{\sigma}\,\{\oz,\oz\}_{\sigma}\,\op \nn \\
S_{\sigma}(\oz)=m_{+}\,\{\oz,\oz\}_{\sigma}\,\om\wedge\op, &\qquad\qquad S_{\sigma}(\om\wedge\op)=-2m_{+}(q^8/\lambda_{(2)}^+)\{\om,\om\}_{\sigma}\,\{\op,\op\}_{\sigma}\,\oz 
\label{Hoq3}
\end{align}
with a normalisation condition $m_+^2\,{\rm det}_{\sigma}\gm\,=\,-{\rm sgn} \gamma$.\bigskip

\bigskip
\bigskip

\item 
Given the quantum tangent space $\mathcal{X}_{\cq_4}$, exact one forms are
\begin{align*}
\dd a=-q\,c^*\omega_{+}\,+\,a\,\omega_{z},& \qquad\qquad \dd c=a^*\omega_{+}\,+\,c\,\omega_{z}, 
\\
\dd a^*=c\,\omega_-\,-\,q\,a^*\omega_{z},& \qquad\qquad \dd c^*=-q^{-1}a\,\omega_{-}\,-\,q\,c^*\omega_z;
\end{align*}
with a braiding:
\begin{align*}
\sigma(\omega_{\pm}\otimes\omega_{\pm})=\omega_{\pm}\otimes\omega_{\pm},& \\
\sigma(\omega_z\otimes\omega_z)=\omega_z\otimes\omega_z\,+\,\frac{1-q}{1+q}\,(\omega_{-}\otimes\omega_+\,-\,\omega_{+}\otimes\omega_{-}),&
\\
\sigma(\omega_+\otimes\omega_-)=q^2\omega_-\otimes\omega_+\,+\,(1-q^2)\omega_+\otimes\omega_-,&\qquad\qquad\sigma(\omega_{-}\otimes\omega_+)=\omega_+\otimes\omega_- \\
\sigma(\oz\ot\om)=q^2\om\ot\oz\,+\,(1-q^2)\oz\ot\om,& \qquad\qquad \sigma(\om\ot\oz)=\oz\ot\om
\\
\sigma(\op\ot\oz)=q^2\oz\ot\op\,+\,(1-q^2)\op\ot\oz,& \qquad\qquad \sigma(\oz\ot\op)=\op\ot\oz.
\end{align*}
The hermitian structure and wedge products read over left-invariant two forms: 
 \begin{align*}
 (\omega_-\wedge\omega_+)^*&=-\om\wedge\op=q^{-2}\omega_{+}\wedge\omega_{-}, \\
 (\omega_-\wedge\omega_z)^*&=-\oz\wedge\op=q^{-2}\omega_{+}\wedge\omega_{z}, \\
(\omega_+\wedge\omega_z)^*&=-\oz\wedge\om=q^{2}\omega_{-}\wedge\omega_{z};
\end{align*}
so that the volume form $\theta_+$ is again a quantum multiple of the classical one:
\beq 
\label{vol4}
\theta_+=q^2(\omega_-\otimes(\omega_{+}\otimes\omega_{z}\,-
\omega_z\otimes\omega_{+})\,+\,\op\otimes(\omega_z\otimes\omega_{-}\,-\,\om\ot\oz) 
\,+\,\oz\ot(\om\ot\op\,-\,\op\ot\om)),
\eeq
giving the following expression $\gm(\theta_+,\theta_+)=-6q^4(\alpha\beta\gamma).$
It is
$$
\{\om,\om\}_{\sigma}\,=\,-q^4\,\alpha, \qquad
\{\op,\op\}_{\sigma}\,=\,-\beta, \qquad
\{\oz,\oz\}_{\sigma}\,=\,-q^{2}\,\gamma.
$$
The set of  $S_{\sigma^{\pm}}$-real and symmetric contraction is 
\beq 
\gm\,\in\,G_{S}\qquad\Leftrightarrow\qquad\{\alpha\,=\,\beta\,\in\,\IR,\, \gamma\,\in\,\IR\};
\label{coq4}
\eeq 
while the conditions of $T_{\sigma}$-reality and symmetry for the same contraction are
$\gm\,\in\,G_{T}\Leftrightarrow\{\alpha\,=\,q^{-4}\beta\,\in\,\IR,\, \gamma\,\in\,\IR\}$.
The Hodge operator is:
\begin{align}
S_{\sigma}(1)=\mu_+, &\qquad\qquad S_{\sigma}(\mu_{+})=-{\rm sgn}(\gamma), \nn \\
S_{\sigma}(\om)=m_+\,q^{-2}\{\om,\om\}_{\sigma}\,\om\wedge\oz, &\qquad\qquad S_{\sigma}(\om\wedge\oz)=2m_+(q^{-2}/\lambda_{(2)}^+)\{\om,\om\}_{\sigma}\,\{\oz,\oz\}_{\sigma}\, \om \nn \\
S_{\sigma}(\op)=-m_+\,\{\op,\op\}_{\sigma}\,\op\wedge\oz, &\qquad\qquad S_{\sigma}(\op\wedge\oz)=-2m_{+}(q^{4}/\lambda_{(2)}^+)\{\op,\op\}_{\sigma}\,\{\oz,\oz\}_{\sigma}\,\op \nn \\
S_{\sigma}(\oz)=-m_+\,\{\oz,\oz\}_{\sigma}\,\om\wedge\op, &\qquad\qquad S_{\sigma}(\om\wedge\op)=-2m_+(1/\lambda_{(2)}^+)\{\om,\om\}_{\sigma}\,\{\op,\op\}_{\sigma}\,\oz 
\label{Hoq4}
\end{align}
with a normalisation condition $m_+^2\,{\rm det}_{\sigma}\gm\,=\,-{\rm sgn} \gamma$. 
\bigskip
\bigskip
\item 
Again we refer to what already noticed, and do not explicitly report the results concerning the calculus generated by $\cq_{5}$ since they can be obtained by the those of the previous calculus by mapping $q\,\to\,-q$.
Once more, the relations \eqref{Hoq4} (and then the \eqref{coq4}) being invariant by this reflection, the Hodge duality we obtain is exactly the previous one.  
\bigskip
\bigskip

\item
Given the quantum tangent space $\mathcal{X}_{\cq_6}$ the exact one forms are
\begin{align*}
\dd a=-q\,c^*\omega_{+}\,-\,a\,\omega_{z},& \qquad\qquad \dd c=a^*\omega_{+}\,-\,c\,\omega_{z}, 
\\
\dd a^*=c\,\omega_-\,+\,q^{4}a^*\omega_{z},& \qquad\qquad \dd c^*=-q^{-1}a\,\omega_{-}\,+\,q^{4}c^*\omega_z;
\end{align*}
and the braiding reads
\begin{align*}
\sigma(\omega_{a}\otimes\omega_{a})=\omega_{a}\otimes\omega_{a},& \qquad\qquad{\rm for}\,a=\pm,z
\\
\sigma(\omega_-\otimes\omega_+)=q^{-4}\omega_+\otimes\omega_-\,+\,(1-q^2)\omega_-\otimes\omega_+\,+\,q^2(q^2-1)\oz\ot\oz,&
\\ \sigma(\omega_{+}\otimes\omega_-)=q^6\omega_-\otimes\omega_+\,-\,q^6(q^2-1)\oz\ot\oz  &\\
\sigma(\om\ot\oz)=q^{-2}\oz\ot\om\,+\,(1-q^2)\om\ot\oz,& \qquad\qquad \sigma(\oz\ot\om)=q^4\om\ot\oz
\\
\sigma(\oz\ot\op)=q^{-2}\op\ot\oz\,+\,(1-q^2)\oz\ot\op,& \qquad\qquad \sigma(\op\ot\oz)=q^4\oz\ot\op.
\end{align*}
The hermitian structure over left-invariant two forms and the wedge product antisymmetry are
 \begin{align*}
 (\omega_-\wedge\omega_+)^*&=-\omega_-\wedge\omega_{+}=q^{-4}\omega_{+}\wedge\omega_{-}, \\
 (\omega_-\wedge\omega_z)^*&=-\omega_z\wedge\omega_{+}=q^{-2}\omega_{+}\wedge\omega_{z}, \\
(\omega_+\wedge\omega_z)^*&=-\omega_z\wedge\omega_{-}=q^{2}\omega_{-}\wedge\omega_{z};
\end{align*}
and the volume form 
\begin{align}
&\theta_+=q^2\omega_-\otimes(q^2\omega_{+}\otimes\omega_{z}\,-q^6
\omega_z\otimes\omega_{+})\nn \\ &\qquad\,+\,q^{-6}\op\otimes(\omega_z\otimes\omega_{-}\,-\,q^4\om\ot\oz) 
\,+\,q^4\oz\ot(\om\ot\op\,-\,q^{-6}\op\ot\om\,+\,(1-q^2)\oz\ot\oz),
\label{vol6}
\end{align}
with  $\gm(\theta_+,\theta_+)\,=\,-\,q^2\gamma (6\,\alpha\,\beta-(1-q^2)^2\gamma^2).$
For the scalar product \eqref{sdpm} one has
$$
\{\om,\om\}_{\sigma}\,=\,-\alpha, \qquad
\{\op,\op\}_{\sigma}\,=\,-q^4\,\beta, \qquad
\{\oz,\oz\}_{\sigma}\,=\,-q^{2}\,\gamma.
$$
The conditions of reality and symmetry of the  contraction \eqref{g1fo} with respect to $S_{\sigma}$ can be expressed by
\beq 
\gm\,\in\,G_{S}\qquad\Leftrightarrow\qquad\{\alpha\,=\,-\ii\,q^6\,\xi, \,\beta\,=\,\ii\,q^4\rho, \, (q^2-1)\gamma\,=\,\pm2\,q^{-2}\xi; \,0\,\neq\,\xi\,\in\,\IR\},
\label{coq6}
\eeq
while 
$\gm\,\in\,G_{T}\Leftrightarrow\{\alpha\,=\,\ii\,\rho, \,\beta\,=\,-\ii\,q^6\rho, \, (q^2-1)\gamma\,=\,\pm2\,\rho; \,0\,\neq\,\rho\,\in\,\IR\}$.
Since these conditions appear counterintuitive, we report the expression that the determinant of the contraction for symmetric and real contractions acquires, namely
\beq
\label{sdd6}
{\rm det}_{\sigma^{\pm}}\gm=\,-\,\frac{2q^8}{\lambda_{(3)}^{\pm}}\,\gamma\,\rho^2.
\eeq
The corresponding Hodge operator is:
\begin{align}
S_{\sigma}(1)=\mu_+, &\qquad\qquad S_{\sigma}(\mu_+)=-{\rm sgn}(\gamma), \nn \\
S_{\sigma}(\om)=m_+\,q^{-4}\{\om,\om\}_{\sigma}\,\om\wedge\oz, &\qquad\qquad S_{\sigma}(\om\wedge\oz)=2m_+(q^{-4}/\lambda_{(2)}^+)\{\om,\om\}_{\sigma}\,\{\oz,\oz\}_{\sigma}\,\om \nn \\
S_{\sigma}(\op)=-m_+\,\{\op,\op\}_{\sigma}\,\op\wedge\oz, &\qquad\qquad S_{\sigma}(\op\wedge\oz)=-2m_+(q^4/\lambda_{(2)}^+)\{\op,\op\}_{\sigma}\,\{\oz,\oz\}_{\sigma}\,\op \nn \\
S_{\sigma}(\oz)=-m_+\,\{\oz,\oz\}_{\sigma}\,\om\wedge\op, &\qquad\qquad \nn \\
 S_{\sigma}(\om\wedge\op)&=-(m_+/\lambda_{(2)}^+)(2q^{-2}\{\om,\om\}_{\sigma}\{\op,\op\}_{\sigma}-(q(q^2-1)(\{\oz,\oz\}_{\sigma})^2)\oz 
\label{Hoq6}
\end{align}
with a normalisation condition $m_+^2\,{\rm det}_{\sigma}\gm\,=\,-{\rm sgn}\, \gamma$. 
Due to the symmetry and reality conditions \eqref{coq6}, the Laplacians  \eqref{laeq}  have for this calculus a spectrum which is not real. This is a condition that characterises this calculus and the corresponding symmetric and real tensors $\gm$ as in \eqref{metr} with respect to the others.

\end{enumerate}

\noindent Once we have the explicit families of  real and symmetric tensors $\gm$ in \eqref{metr} for any of the calculi in $\ck$,  we further address a question arose in the previous pages, namely we wonder whether it is possible  to set a condition eventually selecting between the notions of $T_{\sigma^{\pm}}$ and $S_{\sigma^{\pm}}$ symmetry on one side, and even  among the calculi considered above on the other.

Our approach is straightforward.  
We  consider the set $\mathfrak{G}$ of rank 2 tensors \eqref{metr} satisfying the  symmetry condition $\gm\,=\,\gm\,\circ\,\sigma$ (the meaning of this condition has been extensively analysed when the braiding is the canonical one associated to a bicovariant calculus), and we compare it with $G_{S}, \,G_{T}$. A direct inspection shows that $\mathfrak{G}\,\neq\,G_T$ for any calculus in $\ck$, while $\mathfrak{G}\,=\,G_S$ \emph{only for} the calculi (1,2) and (4,5) following the list above (i.e. those defined in  \eqref{q1}, \eqref{q2}, \eqref{q4}, \eqref{q5}). We call this subset $\tilde{K}\,\subset\,\ck$. 
It is interesting to notice that the Woronowicz calculus does not fulfill this condition; the only calculi satisfying this condition are those, whose volume form is a multiple of the classical one, as it can be immediately seen from the explicit expressions in \eqref{volw}, \eqref{vol1}, \eqref{vol3}, \eqref{vol4}, \eqref{vol6}. 

We shall explore these calculi from a possible different perspective  by studying an extension to the homogeneous space $\sq$ of the formalism giving Hodge dualities on  $\SU$.

\section{Hodge operators over the standard Podle\'s sphere}
\label{shosf}

In section \ref{ss:sesu} we introduced the standard Podle\'s sphere $\Asq$ as the subalgebra of $\U(1)$-coinvariant elements in $\ASU$ by the  coaction given in \eqref{cancoa}. 
As a set of generators for the algebra $\Asq$ we consider
\beq
\label{gps}
B_{-}=-ac^*, \qquad B_{+}=qca^*, \qquad B_{0}=\frac{q^2}{1+q^2}-q^2cc^*,
\eeq
with  $B_{0}^*=B_{0},\,B_{+}^*=-qB_{-}$\footnote{They 
 satisfy the algebraic relations:
\begin{align}
(1+q^{-2})(B_{-}B_{+}+q^2B_{+}B_{-})&=q((1+q^{-2})^2B_{0}^2-1), \nn \\
q(B_{-}B_{+}-B_{+}B_{-})+(q^{-2}-q^2)B_{0}^2&=(1-q^2)B_{0}, \nn \\
(1+q^{-2})(B_{-}B_{0}-q^2B_{0}B_{-})&=(1-q^2)B_{-}, \nn \\
(1+q^{-2})(B_{0}B_{+}-q^2B_{+}B_{0})&=(1-q^2)B_{+}
\label{alsf}
\end{align}
The isomorphism (compatible with the $*$ anti-hermitian conjugation) to the  algebra generated by $\{e_{\pm1}, e_{0}\}$ with relations given in (1)-(4) from \cite{as94} (for the real form $\SU$ of $\mathrm{SL}_{q}(2)$ it is $e_{+1}^*=e_{-1}$) is given by: 
\begin{align}
(1+q^{-2})B_{0}&\mapsto\,e_{0} \nn \\
B_{-}&\mapsto\pm i\,e_{+1} \nn \\
B_{+}&\mapsto\pm iq\,e_{-1} 
\label{iss}
\end{align}
with the identification 
\beq 
\lambda=(1-q^2), \qquad \rho=1.
\label{lrh}
\eeq
}.

Apply now the formalism developed in \cite{bm}: given the 3d calculi in $\ck$ characterised by $\mathcal{Q}_{a}\,\subset\ker\,\varepsilon_{\SU}$, the position  $\pi(\mathcal{Q}_{a})\,\subset\,\ker\,\varepsilon_{\U(1)}$ defines -- from \eqref{qprp} -- a calculus on $\U(1)$. For $a\,\neq\,6$ this calculus on $\U(1)$ is bicovariant and 1 dimensional -- we call \emph{projectable} such 3d calculi on $\SU$, and denote them by $\ck_{\pi}\,\subset\,\ck$ -- and notice that  the  position $\pi(\mathcal{Q}_{6})$ would induce on $\U(1)$  a trivial (0-dimensional) calculus. 
The restriction to $\sq$ of the projectable calculi gives left-covariant calculi: in terms of the generators \eqref{gps} the first order  part of such calculi (that is, for any of the realizations of the operator $\dd$) is  characterized as the left covariant bimodule given by the quotient  $\Gamma(\sq)\,=\,\Asq\,\{\dd B_{\pm},\,\dd B_{0}\}/\omega_0$ where one has defined 
\beq 
\label{doz}
\omega_0\,=\,q^{-1}B_-\dd B_+\,+\,q B_{+}\dd B_-\,-(1+q^{-2})B_0\dd B_0,
\eeq
while their higher order part is given as the   quotient of  the tensor products $\Gamma^{\otimes}(\sq)$ by the differential ideal with generators  $\{\omega_0,\,\dd\omega_0\}$. Such a characterization allows to understand that 
all these calculi on $\sq$ are  isomorphic to the well known 2d left covariant calculus described by Podle\'s in \cite{po92} (the proof of this equivalence is straightforward, mimicking the one  in \cite[\S 3.4]{schm99} which holds for all the projectable calculi over $\SU$). 

It is moreover immediate to prove that this setting describes the geometry of   $\U(1)$ Hopf fibrations over the quantum  sphere $\sq$ with compatible calculi. 
This compatibility allows to meaningful recover that the quantum tangent space $\mathcal{X}_{\pi(\mathcal{Q}_a)}$ associated to the 1 dimensional calculs over $\U(1)$ is vertical for the fibrations, while the 2d exterior algebra over $\sq$ is given by horizontal and $\U(1)$-coequivariant exterior forms on $\SU$ (more details can be found for example in \cite{buao}).  Adopting the so called frame bundle approach \cite{maj} we  write the exterior algebras over the Podle\'s sphere (corresponding to any of the projectable calculi over $\SU$) as
\beq
\label{fbsq}
\Gamma_{\sigma}({\sq}) = \Asq \oplus \left(\cl_{-2} \omega_{-}
\oplus \cl_{+2} \omega_{+} \right) \oplus \Asq \omega_{-}\wedge\omega_{+} \eeq
(i.e. in a more complete  -- and heavier -- notation one should write $\Gamma^{(a)}_{\sigma}(\sq)$ and $\omega_{\pm}^a$ with $(1,\ldots,7)\,\ni\,a\,\neq\,6$.   We also remark that one has an isomorphism $\Gamma_{\sigma^-}(\sq)\,\sim\,\Gamma_{\sigma}(\sq)$ for any projectable calculus). 
This expression is the counterpart in the quantum setting of the classical    \eqref{hoscl}: the sets $\cl_{\pm2}\omega_{\pm}$ are not free $\Asq$-bimodules, the top form bimodule \emph{does} indeed have a 1 dimensional free $\Asq$ left invariant basis. Given the left invariant volume 2-forms
\beq 
\label{vos2}
\check{\mu}_{+}\,=\,\ii\,\check{m}_+\om\wedge\op\,\qquad\qquad
\check{\mu}_{-}\,=\,\ii\,\check{m}_-\om\vee\op\,\qquad\qquad
\eeq 
with $\check{m}_{\pm}\,\in\,\IR$ so to have $\check{\mu}_{\pm}^*\,=\,\check{\mu}_{\pm}$, we define the operators $\check{S}_{\pm}\,:\,\Gamma_{\sigma^{\pm}}^{k}(\sq)\,\to\,\Gamma^{2-k}_{\sigma^{\pm}}(\sq), \, (k=0,1,2)$ via a tensor $\gm\,\in\,G_{S}$ (i.e. $\gm$ is  $S$-symmetric and real on $\SU$) in terms of the decomposition \eqref{fbsq},
\beq
\label{defSs}
\check{S}_{\sigma^{\pm}}(\omega)\,=\,\frac{1}{\lambda^{\mp}_{(k)}}\,\gm(\omega, \check{\mu}_{\pm}).
\eeq
This contraction operator clearly retains the left $\ASU$-linearity in the l.h.s. factor.  The scale is fixed by the natural normalization condition:
\beq
\label{norma}
\check{S}_{\sigma^{\pm}}^2(1)\,=\,{\rm sgn}\,\{\,\gm(\ii\,\om\wedge\op,\ii\,\om\wedge\op)\}.
\eeq
Notice that only for the projectable calculi the braiding over $\SU$ consistently restricts to a braiding among 1-forms on $\sq$ as in the classical setting, so that $\check{S}_{\sigma^{\pm}}$ is well defined; for those calculi the spectral decomposition of the antisymmetrisers on $\sq$ is given  by the restriction of the one on $\SU$ with the same eigenvalues  ${\lambda^{\pm}_{(k)}}$\footnote{It is moreover possible to prove that the element $$\Gamma^{\otimes2}(\sq)\,\ni\,\phi\,=\,q^{-1}\dd B_-\otimes\dd B_+\,+\,q \dd B_{+}\otimes \dd B_-\,-(1+q^{-2})\dd B_0\otimes \dd B_0$$ is the only (up to scalars) to be left invariant and symmetric for any of the resulting braidings: this result extends then the result proved in the  proposition 4.2 in \cite{maj-ri-02}.  }. 

Since for any choice of a projectable calculus on $\SU$ we have a specific realisation  of the unique  2d left covariant $\Gamma(\sq)$,   
what we introduce in \eqref{defSs} are  different (because they come from different calculi on $\SU$ and corresponding different $G_{S}$ sets of symmetric and real tensors $\gm$) contraction operators acting on the same exterior algebra. 

Which is their degeneracy? The only non trivial behaviour to analyze is how they act on 1-forms: a direct proof shows that $(\check{S}_{\sigma^{\pm}})^2$ is a multiple of the identity operator on 1-forms on $\sq$ \emph{only} for those calculi on $\SU$ for which one has that $G_{S}\,=\,\mathfrak{G}$, namely the calculi in $\tilde{K}\,\subset\,\ck$. 
Following the same approach we used for $\SU$, we  shall then define  those contraction operators  \eqref{defSs} corresponding to the calculi in $\tilde{\ck}$ on $\SU$ as the Hodge operators on the unique left covariant 2d calculus on the Podle\'s sphere. 

We explicitly write  these two pairs of Hodge operators.
\begin{itemize}
\item We realize  the 2d differential calculus \eqref{fbsq} on $\sq$ as a restriction of the calculus \eqref{q1} on $\SU$, characterized by $\mathcal{Q}_{1}\,\subset\,\ker\,\varepsilon_{\SU}$.  
The action of the Hodge operator turns out to be (we know from \eqref{coq1} that it depends on the real parameter $\alpha\,\neq\,0$):
\begin{align}
\check{S}_{\sigma}(1)&=\,\ii\check{m}_{+}\,\om\wedge\op, \nn \\
\check{S}_{\sigma}(x_{\pm}\,\omega_{\pm})&=\,\pm\ii\check{m}_{+}\,q^2\,\alpha\,x_{\pm}\,\omega_{\pm},\nn\\
\check{S}_{\sigma}(\om\wedge\op)&=\,-\ii\check{m}_{+}\,2\,q^4\,\alpha^2/\lambda_{(2)}^-,
\label{ss1z}
\end{align}
where $x_{\pm}\,\in\,\cl_{\pm2}$, while the normalization reads $\check{m}^2_{+}\,2\,q^4\,\alpha^2/\lambda_{(2)}^-\,=\,1$. Given the isomorphism $\Gamma_{\sigma}(\sq)\,\sim\,\Gamma_{\sigma^-}(\sq)$, from \eqref{isospec} it is immediate to see that
\begin{align} 
\check{S}_{\sigma^-}(1)&=\,\ii q^{-2}\check{m}_{-}\,\om\wedge\op, \nn \\
\check{S}_{\sigma}(x_{\pm}\,\omega_{\pm})&=\,\pm\ii\check{m}_{-}\,\alpha\,x_{\pm}\,\omega_{\pm},\nn\\
\check{S}_{\sigma}(\om\wedge\op)&=\,-\ii\check{m}_{-}\,2\,q^2\,\alpha^2/\lambda_{(2)}^+,
\label{ss1mz}
\end{align}
with $\check{m}^2_{-}\,2\,\alpha^2/\lambda_{(2)}^+\,=\,1$. We recall that the Hodge operators on $\sq$ corresponding to those  defined on $\SU$ via the calculus (2) (that is in \eqref{q2}) coincide with those above, since it is obtained by mapping $q\,\rightarrow\,-q$. 

\bigskip
\item
In the same way, if we realize   the 2d differential calculus \eqref{fbsq} on $\sq$ as a restriction of the calculus \eqref{q4} on $\SU$, characterized by $\mathcal{Q}_{4}\,\subset\,\ker\,\varepsilon_{\SU}$, then the Hodge operator we obtain is (recall  from \eqref{coq4} that it also depends on the real parameter $\alpha\,\neq\,0$):
\begin{align}
\check{S}_{\sigma}(1)&=\,\ii\check{m}_{+}\,\om\wedge\op, \nn \\
\check{S}_{\sigma}(x_{\pm}\,\omega_{\pm})&=\,\pm\ii\check{m}_{+}\,\alpha\,x_{\pm}\,\omega_{\pm},\nn\\
\check{S}_{\sigma}(\om\wedge\op)&=\,-\ii\check{m}_{+}\,2\,\alpha^2/\lambda_{(2)}^-,
\label{ss4z}
\end{align}
where $x_{\pm}\,\in\,\cl_{\pm2}$, while the normalization reads $\check{m}^2_{+}\,2\,\alpha^2/\lambda_{(2)}^-\,=\,1$. Given the isomorphism $\Gamma_{\sigma}(\sq)\,\sim\,\Gamma_{\sigma^-}(\sq)$, from \eqref{isospec} it is immediate to see that
\begin{align} 
\check{S}_{\sigma^-}(1)&=\,\ii q^{-2}\check{m}_{-}\,\om\wedge\op, \nn \\
\check{S}_{\sigma}(x_{\pm}\,\omega_{\pm})&=\,\pm\ii q^{-2}\check{m}_{-}\,\alpha\,x_{\pm}\,\omega_{\pm},\nn\\
\check{S}_{\sigma}(\om\wedge\op)&=\,-\ii\check{m}_{-}\,2\,q^{-2}\,\alpha^2/\lambda_{(2)}^+,
\label{ss1mz}
\end{align}
with $\check{m}^2_{-}\,2q^{-4}\,\alpha^2/\lambda_{(2)}^+\,=\,1$. Once more we have that the Hodge operators on $\sq$ corresponding to those  defined on $\SU$ via the calculus (5) (i.e. \eqref{q5}) coincide with  the ones above.

\end{itemize}

We close this analysis by describing the Laplacian operators introduced on $\sq$ by the Hodge operators above. A direct calculation shows that one has  (due to the normalization condition given in \eqref{norma})
\beq
\label{lasq2}
\Box_{\sq}\,f\,=\,q\,(EF\,+\,FE)\,\lt f
\eeq
(on $f\,\in\,\Asq$) for any of the calculi in $\tilde{\ck}$  and for any realization of the exterior algebra in terms of the braidings $\sigma^{\pm}$.  The action of this operator moreover coincides with the action of the Laplacian operator on $\sq$ obtained in \cite{ale10} following a different formulation.

\bigskip
\bigskip
\bigskip
\bigskip

{\bf Acknowledgments.}  This paper underwent various revisions during the last months. For any new version of it I am indebted to friends and colleagues for their precious feedback. Giovanni Landi has been along this research a wonderful guidance; Francesco D'Andrea, Istvan Heckenberger, Debashish Goswami gifted me some of their insights on the subjects; Sergio Albeverio, Francesco Bonechi, Yuri I. Manin, Giuseppe Marmo, Gianluca Panati, Sylvie Paycha, Alessandro Teta,  gave me the opportunity to present part of this paper; Gianfausto Dell'Antonio and Detlef D\"urr mentored and supported me during the last year spent at L.M.U. in M\"unchen. I express them all my deep  gratitude. 
It is a pleasure for me to acknowledge a financial support of the H.C.M in Bonn.

\end{document}